\documentclass[reqno]{amsart}						
\usepackage{microtype}
\usepackage{fullpage}

\usepackage{amsmath}						
\usepackage{amssymb}						
\usepackage{amsfonts}						
\usepackage{amsthm}	
\usepackage[foot]{amsaddr}						
\usepackage{array}

\usepackage{mathtools}						
\mathtoolsset{%
}

\usepackage[utf8]{inputenc}							
\usepackage[T1]{fontenc}						

\usepackage[proportional,lining]{libertine}
\usepackage[libertine,libaltvw,cmintegrals,varbb]{newtxmath}

\usepackage{dsfont}							

\usepackage[
cal=cm,
]
{mathalfa}

\usepackage[labelfont={bf,small},labelsep=colon,font=small]{caption}	
\usepackage[dvipsnames,svgnames]{xcolor}						
\colorlet{MyBlue}{DodgerBlue!75!Black}
\colorlet{MyGreen}{DarkGreen!85!Black}

\usepackage{subcaption}
\usepackage{tikz}							
\usetikzlibrary{calc,patterns,matrix}
\usepackage{algorithmic}
\usepackage{algorithm}

\newcommand\overmat[2]{%
  \makebox[0pt][l]{$\smash{\color{black}\overbrace{\phantom{%
    \begin{matrix}#2\end{matrix}}}^{\text{\color{black}#1}}}$}#2}

\usepackage{tikz}
\usepackage{pgfplots}
\usepackage{mathrsfs}

\usetikzlibrary{arrows}

\usetikzlibrary{plotmarks,patterns}
\usetikzlibrary{matrix}

\usepackage{acronym}						
\usepackage{booktabs}						
\usepackage{latexsym}						
\usepackage{paralist}						
\usepackage{xspace}							
\usepackage{multirow}

\usepackage{natbib}
 \bibpunct[, ]{[}{]}{,}{n}{}{,}%

\usepackage{hyperref}
\hypersetup{
final,
colorlinks=true,
linktocpage=true,
pdfstartview=FitH,
breaklinks=true,
pdfpagemode=UseNone,
pageanchor=true,
pdfpagemode=UseOutlines,
plainpages=false,
bookmarksnumbered,
bookmarksopen=false,
bookmarksopenlevel=1,
hypertexnames=true,
pdfhighlight=/O,
urlcolor=Maroon,linkcolor=MyBlue!60!black,citecolor=DarkGreen!70!black,	
pdftitle={},
pdfauthor={},
pdfsubject={},
pdfkeywords={},
pdfcreator={pdfLaTeX},
pdfproducer={LaTeX with hyperref}
}

\def\EMAIL#1{\email{\href{mailto:#1}{\texttt{\upshape #1}}}}

\theoremstyle{plain}
\newtheorem{theorem}{Theorem}						
\newtheorem*{theorem*}{Theorem}						
\newtheorem*{corollary*}{Corollary}						
\newtheorem{lemma}{Lemma}					

\theoremstyle{definition}
\newtheorem{definition}{Definition}				
\newtheorem*{definition*}{Definition}					
\newtheorem{assumption}{Assumption}					

\theoremstyle{remark}
\newtheorem{remark}{Remark}						
\newtheorem*{remark*}{Remark}						
\newtheorem{example}{Example}						
\newtheorem*{example*}{Example}						

\DeclareMathOperator{\diag}{diag}

\newcommand{\algo}{\texttt{RPSD}\xspace}
\newcommand{\pgd}{\texttt{PGD}\xspace}
\newcommand{\rpcd}{\texttt{RPCD}\xspace}
\newcommand{\adaalgo}{\texttt{ARPSD}\xspace}

\newcommand{\sega}{\texttt{SEGA}\xspace}

\DeclareMathOperator*{\argmin}{argmin}
\newcommand{\prox}{\mathbf{prox}}

\newcommand{\rank}{\mathrm{rank}}

\newcommand{\EE}{\mathbb{E}}

\newcommand{\RR}{\mathbb{R}}
\newcommand{\Prob}{\mathbb{P}}
\newcommand{\C}{\mathcal{C}}

\newcommand{\M}{\mathcal{M}}
\renewcommand{\S}{\mathsf{S}_\M}
\newcommand{\SC}{\mathsf{S}_{\NC}}
\newcommand{\NC}{\mathcal{M}}

\newcommand{\Sel}{\mathfrak{S}}

\newcommand{\bP}{{\mathsf{P}}}
\newcommand{\bQ}{{\mathsf{Q}}}

\newcommand{\FF}{\mathcal{F}}
\newcommand{\Up}{\mathsf{L}}

\begin{document}

\title{Proximal Gradient methods with Adaptive Subspace Sampling}


\author[D.~Grishchenko]{Dmitry Grishchenko$^{\star}$}

\author[F.~Iutzeler]{Franck Iutzeler$^{\star}$}

\address{$^{\star}$\,%
Univ. Grenoble Alpes, Laboratoire Jean Kuntzmann}

\author[J.~Malick]{Jérôme Malick$^{\circ}$}
\address{$^{\circ}$\,%
CNRS, Laboratoire Jean Kuntzmann}

\EMAIL{firstname.lastname@univ-grenoble-alpes.fr}

\maketitle

\begin{abstract}
    Many applications in machine learning or signal processing involve nonsmooth optimization problems. This nonsmoothness brings a low-dimensional structure to the optimal solutions. In this paper, we propose a randomized proximal gradient method harnessing this underlying structure. We introduce two key components: i) a random subspace proximal gradient algorithm; ii) an identification-based sampling of the subspaces. Their interplay brings a significant performance improvement on typical learning problems in terms of dimensions explored.
\end{abstract}

\section{Introduction}

In this paper, we consider composite optimization problems of the form
\begin{equation}\label{eq:main_problem}
	\min_{x\in \RR^n}~f(x)+g(x)
\end{equation}
where $f$ is convex and differentiable, and $g$ is convex and nonsmooth. 
This type of problem appears extensively in signal processing and machine learning applications; we refer to e.g.\;\cite{candes2008enhancing}, \cite{combettes2011proximal}, \cite{bach2012optimization}, among a vast literature. Large scale applications in these fields call for first-order optimization, such as proximal gradient methods
(see e.g.\;the recent\;\cite{teboulle2018simplified}) and coordinate descent algorithms (see e.g.\;the review\;\cite{wright2015coordinate}).

In these methods, the use of a proximity operator to handle the nonsmooth part $g$ plays a prominent role, as it typically enforces some ``sparsity'' structure on the iterates and eventually on optimal solutions, see e.g.\;\cite{vaiter-model-linear15}. For instance, the popular $\ell_1$-norm regularization ($g=\|\cdot\|_1$) promotes optimal solutions with a few nonzero elements, and its associated proximity operator (called soft-thresholding, see\;\cite{donoho1995noising}) zeroes entries along the iterations. This is an example of \emph{identification}: in general, the iterates produced by proximal algorithms eventually reach some sparsity pattern close to the one of the optimal solution. For $\ell_1$-norm regularization, this means that after a finite but unknown number of iterations the algorithm ``identifies'' the final set of non-zero variables. This active-set identification property is typical of constrained convex optimization (see e.g.\;\cite{wright1993identifiable}) and nonsmooth optimization (see e.g.\;\cite{hare2004identifying}).

The study of identification dates back at least to \cite{bertsekas1976goldstein} who showed that the projected gradient method identifies a sparsity pattern when using non-negative constraints. Such identification has been extensively studied in more general settings; we refer to \cite{burke1988identification}, \cite{lewis2002active}, \cite{drusvyatskiy2013optimality} or the recent \cite{lewis2018partial}, among other references.
Recent works on this topic include: i) extended identification for a class of functions showing strong primal-dual structure, including TV-regularization and nuclear norm \cite{fadili2018sensitivity}; ii) identification properties of various randomized algorithms, such as coordinate descent \cite{wright2012accelerated} and stochastic methods \cite{pmlr-v80-poon18a,fadili2018model,sun2019we}.

The knowledge of the optimal substructure would allow to reduce the optimization problem in this substructure and solve a lower dimension problem. While identification can be guaranteed in special cases  (e.g.\;using duality\;for $\ell_1$-regularized least-squares \cite{ogawa2013safe, fercoq2015mind}), it is usually unknown beforehand and proximal algorithms can be exploited to obtain approximations of this substructure.
After some substructure identification, one could switch to a more sophisticated method, e.g.\,updating parameters of first-order methods (\cite{2015-Liang-InterialFB}). Again, since the final identification moment is not known, numerically exploiting identification to accelerate the convergence of first-order methods has to be done with great care.


In this paper, we propose randomized proximal algorithms leveraging on structure identification: our idea is to sample the variable space according to the structure of $g$. To do so, we first introduce a randomized descent algorithm going beyond separable nonsmoothness and associated coordinate descent methods:
we consider ``subspace descent" extending ``coordinate descent" to generic subspaces.
Then, we use a standard identification property of proximal methods to adapt our sampling of the subspaces with the identified structure. This results in a structure-adapted randomized method with automatic dimension reduction, which 
performs better in terms of dimensions explored compared standard proximal methods and the non-adaptive version.

Though our main concern is the handling of non-separable nonsmooth functions $g$, we mention that our identification-based adaptive approach is different from existing adaptation strategies restricted to the particular case of coordinate descent methods. Indeed, adapting coordinate selection probabilities is an important topic for coordinate descent methods as both theoretical and practical rates heavily depend on them (see e.g.\;\cite{richtarik2014iteration,necoara2014random}).
Though the optimal theoretical probabilities, named importance sampling, often depend on unknown quantities, these \emph{fixed} probabilities can sometimes be computed and used in practice, see\;\cite{zhao2015stochastic,richtarik2016optimal}.
The use of \emph{adaptive} probabilities is more limited; some heuristics without convergence guarantees can be found in \cite{loshchilov2011adaptive,glasmachers2013accelerated}, and greedy coordinates selection are usually expensive to compute \cite{dhillon2011nearest,nutini2015coordinate,nutini2017let}. Bridging the gap between greedy and fixed importance sampling, \cite{perekrestenko2017faster,namkoong2017adaptive,stich2017safe} propose adaptive coordinate descent methods based on the coordinate-wise Lipschitz constants and current values of the gradient. 
The methods proposed in the present paper, even when specialized in the coordinate descent case, are the first ones where the \emph{iterate structure enforced by a non-smooth regularizer} is used to adapt the selection probabilities.

The paper is organized as follows. In Section~\ref{sec:sub_des}, we introduce the formalism for subspace descent methods. First, we formalize how to sample subspaces and introduce a first random subspace proximal gradient algorithm. Then, we show its convergence and derive its linear rate in the strongly convex case. Along the way, we make connections and comparisons with the literature on coordinate descent and sketching methods, notably in the special cases of $\ell_1$ and total variation regularization. In Section~\ref{sec:ada_sub_des}, we present our identification-based adaptive algorithm. We begin by showing the convergence of an adaptive generalization of our former algorithm; next, we show that this algorithm enjoys some identification property and give practical methods to adapt the sampling, based on generated iterates, leading to refined rates. Finally, in Section~\ref{sec:num}, we report numerical experiments on popular learning problems to illustrate the merits and reach of the proposed methods.

\section{Randomized subspace descent}\label{sec:sub_des}

The premise of randomized subspace descent consists in repeating two steps: i) randomly selecting some subspace; and ii) updating the iterate over the chosen subspace. Such algorithms thus extend usual coordinate descent to general sampling strategies, which requires algorithmic changes and an associated mathematical analysis.
This section presents a subspace descent algorithm along these lines for solving \eqref{eq:main_problem}. In Section\;\ref{sec:sub}, we introduce our subspace selection procedure. We build on it to introduce, in Section\;\ref{sec:algo}, our first subspace descent algorithm, the convergence of which is analyzed in Section\;\ref{sec:conv}. Finally, we put this algorithm into perspective in Section\;\ref{sec:comparison} by connecting and comparing it to related work.

\subsection{Subspace selection}\label{sec:sub}

We begin by introducing the mathematical objects leading to the subspace selection used in our randomized subspace descent algorithms. Though, in practice, most algorithms rely on projection matrices, our presentation highlights intrinsic subspaces associated to these matrices; this opens the way to a finer analysis, especially in Section~\ref{sec:ada_algo} when working with adaptive subspaces. 

We consider a family $\C = \{ \C_i \}_i$ of (linear) subspaces of $\mathbb{R}^n$. Intuitively, this set represents the directions that will be \emph{favored} by the random descent; in order to reach a global optimum, we naturally assume that the sum\footnote{In the definition and the following, we use the natural set addition (sometimes called the Minkowski sum): for any two sets $\mathcal{C},\mathcal{D}\subseteq\mathbb{R}^n$, the set  $\mathcal{C} + \mathcal{D}$ is defined as $\{x+y : x\in\mathcal{C},y\in\mathcal{D}\}\subseteq\mathbb{R}^n$.} of the subspaces in a family matches the whole space.

\begin{definition}[Covering family of subspaces]\label{def:cov}
Let $\C = \{ \C_i \}_i$ be a family of subspaces of $\mathbb{R}^n$. We say that $\C$ is \emph{covering} if it spans the whole space, i.e. if $\sum_i \C_i = \mathbb{R}^n$.
\end{definition}

\begin{example}\label{ex:axes}
The family of the axes $\C_i = \{x\in\mathbb{R}^n  : x_j = 0 ~~ \forall j \neq i \}$ for $i=1,..,n$ is a canonical covering family for $\RR^n$.
\end{example}

From a covering family $\C$, we call \emph{selection} the random subspace obtained by randomly choosing some subspaces in $\C$ and summing them. We call \emph{admissible} the selections that include all directions with some positive probability; or, equivalently, the selections to which no non-zero element of $\mathbb{R}^n$ is orthogonal with probability one.

\begin{definition}[Admissible selection]
Let $\C$ be a covering family of subspaces of $\mathbb{R}^n$. A selection $\Sel$ is defined from the set of all subsets of $\C$ to the set of the subspaces of $\mathbb{R}^n$ as
 $$ \Sel(\omega) = \sum_{j = 1}^s \C_{i_j} \qquad\text{ for } \omega = \{ \C_{i_1} ,\ldots, \C_{i_s} \} .$$
 The selection $\Sel$ is \emph{admissible} if  $\Prob[x\in \Sel^\perp]<1$ for all $x\in\mathbb{R}^n\setminus\{0\}$.
\end{definition}

Admissibility of selections appears on spectral properties of the average projection matrix onto the selected subspaces. For a subspace $F \subseteq \mathbb{R}^n$, we denote by $P_F \in\mathbb{R}^{n\times n}$ the orthogonal projection matrix onto $F$. The following lemma shows that the average projection associated with an admissible selection is positive definite; this matrix and its extreme eigenvalues will play a major role in our developments.

\begin{lemma}[Average projection]\label{lm:eligible}
If a selection $\Sel$ is admissible then 
    \begin{equation}\label{eq:pbar}
        \bP := \EE[P_{\Sel}]  \qquad\text{is a positive definite matrix.}
    \end{equation}
In this case, we denote by $\lambda_{\min}(\bP)>0$ and $\lambda_{\max}(\bP)\leq1$ its minimal and maximal eigenvalues.
\end{lemma}

\proof{Proof.}
Note first that for almost all $\omega$, the orthogonal projection $P_{\Sel(\omega)}$ is positive semi-definite, and therefore so is $\bP$. Now, let us prove that if $\bP$ is not positive definite, then $\Sel$ is not admissible. 
Take a nonzero $x$ in the kernel of $\bP$, then
\begin{equation*}
    x^\top\bP x = 0 \iff x^\top \EE[P_{\Sel}]x = 0 \iff   \EE[x^\top P_{\Sel} x] = 0.
\end{equation*}
Since $x^\top P_{\Sel(\omega)}x\geq 0$ for almost all $\omega$, the above property is further equivalent for almost all $\omega$ to 
\begin{equation*}
    x^\top P_{\Sel(\omega)} x = 0 \iff P_{\Sel(\omega)} x = 0 
    \iff   x\in\Sel(\omega)^{\perp} . 
\end{equation*}
Since $x\neq0$, this yields that $x\in\Sel(\omega)^\perp$ for almost all $\omega$ which is in contradiction with $\Sel$ being admissible. Thus, if a selection $\Sel$ is admissible, $ \bP := \EE[P_{\Sel}]$ is positive definite (so $\lambda_{\min}(\bP)>0$). 

Finally, using Jensen's inequality and the fact that $P_{\Sel}$ is a projection, we get $\| \bP x \| = \| \EE[P_{\Sel}] x \|  \leq \EE \| P_{\Sel} x \| \leq \|x\|$, which implies that $\lambda_{\max}(\bP)\leq1$. 
 \endproof

Although the framework, methods, and results presented in this paper allow for infinite subspace families (as in sketching algorithms); the most direct applications of our results only call for finite families for which the notion of admissibility can be made simpler.

\begin{remark}[Finite Subspace Families]
\label{rem:selection}
For a covering family of subspaces $\C$ with a finite number of elements, the admissibility condition can be simplified to  $\Prob[ \C_i \subset \Sel]>0$ for all $i$.

Indeed, take $x\in\mathbb{R}^n\setminus\{0\}$; then, since  $\C$ is covering and $x\neq 0$, there is a subspace $\C_i$ such that $P_{\C_i}x \neq 0$. Observe now that $\C_i \subset \Sel$ yields $P_{\Sel} x \neq 0$ (since $\Sel^\perp \subset  \C_i^\perp $, the property $P_{\Sel}x=0$ would give $P_{\C_i}x=0$ which is a contradiction with $P_{\C_i}x \neq 0$). Thus, we can write
\begin{align*}
    \Prob[x\in \Sel^\perp] &=  \Prob[P_{\Sel} x = 0] = 1 - \Prob[P_{\Sel} x \neq 0]  \leq 1 - \Prob[ \C_i \subset \Sel ] < 1.
\end{align*}

Building on this property, two natural ways to generate admissible selections from a finite covering family $\C=\{ \C_i \}_{i=1,\ldots,c}$ are:
\begin{itemize}
    \item \emph{Fixed probabilities:} Selecting each subspace $\C_i$ according to the outcome of a Bernoulli variable of parameter $p_i>0$. This gives admissible selections as $\Prob[\C_i \subseteq \Sel] = p_i >0$ for all $i$;
    \item \emph{Fixed sample size:} Drawing $s$ subspaces in $\C$ uniformly at random. This gives admissible selections since $ \Prob[\C_i \subseteq \Sel] = s/c$ for all $i$. 
\end{itemize}
\end{remark}

\begin{example}[Coordinate-wise projections]\label{ex:P}
Consider the family of the axes from Example~\ref{ex:axes} and the selection generated with fixed probabilities as described in Remark~\ref{rem:selection}. The associated projections amount to zeroing entries at random and the average projection $\bP$ is the diagonal matrix with entries $(p_i)$; trivially $\lambda_{\min}(\bP)= \min_i p_i$ and $\leq \lambda_{\max}(\bP) = \max_i p_i$. 
\end{example}

\subsection{A random subspace proximal gradient algorithm} \label{sec:algo}

An iteration of the proximal gradient algorithm 
decomposes in two steps (sometimes called ``forward'' and  ``backward''):
\begin{subequations}\label{eq:pgGradProx}
    \begin{align}
    \label{eq:pgGrad} z^k &= x^{k} - \gamma \nabla f(x^k)\\
        \label{eq:pgProx} x^{k+1} &= \prox_{\gamma g}(z^k)
\end{align}
\end{subequations}
where $\prox_{\gamma g}$ stands for the proximity operator defined as the mapping 
from $\RR^n$ to $\RR^n$
\begin{equation}
    \prox_{\gamma g}(x) = \argmin_{y\in\RR^n}\left\{g(y) + \frac{1}{2\gamma} \|y-x\|_2^2\right\} .
\end{equation}

This operator is well-defined when $g$ is a proper, lower semi-continuous convex function \cite[Def.~12.23]{bauschke2011convex}. Furthermore, it is computationally cheap to compute in several cases, either from a closed form ({e.g.} for $\ell_1$-norm, $\ell_1/\ell_2$-norm, see among others \cite{combettes2007proximal} and references therein), or by an efficient procedure ({e.g.} for the 1D-total variation, projection on the simplex, see \cite{yuan2011efficient,condat2013direct}).

In order to construct a ``subspace'' version of the proximal gradient \eqref{eq:pgGradProx}, one has to determine which variable will be updated along the randomly chosen subspace (which we will call a projected update). Three choices are possible: 
\begin{itemize}
    \item[(a)] a projected update of $x^k$,  i.e. projecting after the proximity operation;
    \item[(b)] a projected update of $\nabla f(x^k)$,  i.e. projecting after the gradient;
    \item[(c)] a projected update of $z^k$,  i.e. projecting after the gradient \emph{step}.
\end{itemize}
Choice (a) has limited interest in the general case where the proximity operator is not separable along subspaces and thus a projected update of $x^k$ still requires the computations of the full gradient. In the favorable case of coordinate projection and $g=\|\cdot\|_1$, it was studied in \cite{qu2016coordinate} using the fact that the projection and the proximity operator commute. Choice (b) is considered recently in \cite{hanzely2018sega} in the slightly different context of sketching. A further discussion on related literature is postponed to Section~\ref{sec:comparison}.

In this paper, we will consider Choice (c), inspired by recent works highlighting that combining iterates usually works well in practice (see \cite{mishchenko2018distributed} and references therein). However, taking gradient steps along random subspaces introduce bias and thus such a direct extension fails in practice. In order to retrieve convergence to the optimal solution of \eqref{eq:main_problem}, we slightly modify the proximal gradient iterations by including a correction featuring the inverse square root of the expected projection denoted by $ \bQ= \bP^{-1/2}$ (note that as soon as the selection is admissible, $ \bQ$ is well defined from Lemma~\ref{lm:eligible}).

Formally, our Random Proximal Subspace Descent algorithm \algo, displayed as Algorithm~\ref{alg:strata_nondis}, replaces \eqref{eq:pgGrad} by 
\begin{equation}\label{eq:grad_step}
    y^k = \bQ\left(x^k - \gamma\nabla f\left(x^k\right)\right)\qquad\text{and}\qquad
    z^{k} = P_{\Sel^k} \left(y^k\right) + (I- P_{\Sel^k} ) \left(z^{k-1}\right).
\end{equation}
That is, we propose to first perform a gradient step followed by a change of basis (by multiplication with the positive definite matrix $\bQ$), giving variable $y^k$; then, variable $z^k$ is updated only in the random subspace $\Sel^k$: to  $P_{\Sel^k} \left(y^k\right)$ in $\Sel^k$, and keeping the same value outside. Note that $y^k$ does not actually have to be computed and only the ``$P_{\Sel^k}\bQ$-sketch'' of the gradient (i.e.\;$ P_{\Sel^k}\bQ\nabla f\left(x^k\right)$) is needed. Finally, the final proximal operation \eqref{eq:pgProx} is performed after getting back to the original space (by multiplication with $\bQ^{-1}$):
\begin{align}\label{eq:prox_step}
   x^{k+1} = \prox_{\gamma g} \left(\bQ^{-1}\left(z^{k}\right)\right).
\end{align}
Contrary to existing coordinate descent methods, our randomized subspace proximal gradient algorithm does not assume that the proximity operator $\prox_{\gamma g} $ is separable with respect to the projection subspaces. Apart from the algorithm of \cite{hanzely2018sega} in a different setting, this is an uncommon but highly desirable feature to tackle general composite optimization problems.  

\begin{algorithm} 
\caption{Randomized Proximal Subspace Descent - \algo}
\label{alg:strata_nondis} 
\begin{algorithmic}[1] 
    \STATE Input:  $\bQ = \bP^{-\frac12}$
    \STATE Initialize $z^0$, $x^1 = \prox_{\gamma g}(\bQ^{-1}(z^0))$
    \FOR{$k=1,\ldots$}
            \STATE $y^k = \bQ\left(x^k - \gamma\nabla f\left(x^k\right)\right)$
            \STATE $z^{k} = P_{\Sel^k} \left(y^k\right) + (I- P_{\Sel^k} ) \left(z^{k-1}\right)$
            \STATE$x^{k+1} = \prox_{\gamma g} \left(\bQ^{-1}\left(z^{k}\right)\right)$
    \ENDFOR
\end{algorithmic}
\end{algorithm}

Let us provide a first example, before moving to the analysis of the algorithm in the next section. 

\begin{example}[Interpretation for smooth problems]
In the case where $g\equiv0$, our algorithm has two interpretations. 
First, using\;$\prox_{\gamma g} = I$, the iterations simplify to 
\begin{align*}
     z^{k+1} &= z^k - \gamma P_{\Sel^k}\bQ\left(\nabla f\left(\bQ^{-1}\left(z^k\right)\right)\right) = z^k - \gamma  P_{\Sel^k}\bQ^2 \underbrace{\bQ^{-1}  \left(\nabla f\left(\bQ^{-1}\left(z^k\right)\right)\right) }_{\nabla f\circ \bQ^{-1} (z^k)}.
\end{align*}
As $\mathbb{E}[P_{\Sel^k}\bQ^2] = I $, this corresponds to a random subspace descent on $ f\circ \left(\bQ^{-1}\right)$ with unbiased gradients. Second, we can write it with the change of variable $u^k =\bQ^{-1} z^k $ as
\begin{align*}
  u^{k+1} &= u^k - \gamma  \bQ^{-1}  P_{\Sel^k} \bQ\left(\nabla f\left(u^k\right)\right).
\end{align*}
As $\mathbb{E} [\bQ^{-1}  P_{\Sel^k}\bQ] = \bP $, this corresponds to random subspace descent on $ f $ but with biased gradient. We note that the recent work \cite{frongillo2015convergence} considers a similar set-up and algorithm; however, the provided convergence result does not lead to the convergence to the optimal solution (due to the use of the special semi-norm). 
\end{example}

\subsection{Analysis and convergence rate}\label{sec:conv}

In this section, we provide a theoretical analysis for \algo, showing linear convergence for strongly convex objectives. Tackling the non-strongly convex case requires extra-technicalities; we thus choose to postpone the corresponding convergence result to the appendix for clarity.

\begin{assumption}[On the optimization problem]\label{hyp:f}
The function $f$ is $L$-smooth and $\mu$-strongly convex and the function $g$ is convex, proper, and lower-semicontinuous.
\end{assumption}

Note that this assumption implies that Problem~\eqref{eq:main_problem} has a unique solution that we denote $x^\star$ in the following.

\begin{assumption}[On the randomness of the algorithm]\label{hyp:main}
Given a covering family $\C=\{ \C_i \}$ of subspaces, we consider a sequence $\Sel^1,\Sel^2,..,\Sel^k$ of admissible selections, which is i.i.d.
\end{assumption}

In the following theorem, we show that the proposed algorithm converges linearly at a rate that only depends on the function properties and on the smallest eigenvalue of $\bP$. We also emphasize that the step size $\gamma$ can be taken in the usual range for the proximal gradient descent.

\begin{theorem}[\algo~convergence rate]\label{th:conv_nondis}
Let Assumptions \ref{hyp:f} and  \ref{hyp:main} hold. Then, for any $\gamma\in(0,2/(\mu+L)]$, the sequence $(x^k)$ of the iterates of \algo converges almost surely to the minimizer  $x^\star$  of \eqref{eq:main_problem} with rate
$$
\EE\left[\|x^{k+1} - x^\star\|_2^2\right] \leq \left(1 - \lambda_{\min}(\bP) \frac{2\gamma\mu L}{\mu + L} \right)^k C,
$$
where $C = \lambda_{\max}(\bP)\|z^0 - \bQ(x^\star - \gamma\nabla f(x^\star))\|_2^2$.
\end{theorem}

To prove this result, we first demonstrate two intermediate lemmas respectively expressing the distance of $z^k$ towards its fixed points (conditionally to the filtration of the past random subspaces $\FF^k = \sigma( \{\Sel_\ell\}_{\ell\leq k} )$), and bounding the increment (with respect to $\| x \|_{\bP}^2 = \langle x , \bP x \rangle $ the norm associated to $\bP$).

\begin{lemma}[Expression of the decrease as a martingale]\label{lm:removing_exp}
From the minimizer  $x^\star$  of \eqref{eq:main_problem}, define the fixed points $z^\star = y^\star = \bQ\left(x^\star - \gamma\nabla f\left(x^\star\right)\right)$ of the sequences $(y^k)$ and $(z^k)$.  If Assumption \ref{hyp:main}~holds, then 
    $$
    \EE\left[\|z^{k} - z^\star\|_2^2\,|\,\FF^{k-1}\right] = \|z^{k-1}-z^\star\|_2^2 + \|y^{k}-y^\star\|_{\bP}^2 -  \|z^{k-1}-z^\star\|_{\bP}^2.
    $$
\end{lemma}

\proof{Proof.}
By taking the expectation on $\Sel^k$ (conditionally to the past), we get 
\begin{align*}
\EE\left[\|z^{k} - z^\star\|_2^2|\,\FF^{k-1}\right]  &= \EE \left[\|z^{k-1} - z^\star + P_{\Sel^k}(y^k-z^{k-1}) \|_2^2 |\,\FF^{k-1}\right]\\
&= \|z^{k-1}-z^\star\|_2^2 + 2 \EE \left[  \langle z^{k-1}-z^\star, P_{\Sel^k}(y^k - z^{k-1} ) \rangle |\,\FF^{k-1}\right] + \EE \left[ \left\| P_{\Sel^k}(y^k-z^{k-1}) \right\|^2 |\,\FF^{k-1}\right] \\
&=\|z^{k-1}-z^\star\|_2^2 + 2\langle z^{k-1}-z^\star, \bP(y^k - z^{k-1} ) \rangle + \EE \left[\langle  P_{\Sel^k}(y^k-z^{k-1}), P_{\Sel^k} (y^k-z^{k-1})\rangle |\,\FF^{k-1}\right]\\
&=\|z^{k-1}-z^\star\|_2^2 + 2\langle z^{k-1}-z^\star, \bP(y^k - z^{k-1} ) \rangle + \EE \left[ \langle y^k-z^{k-1}, P_{\Sel^k} (y^k-z^{k-1})\rangle |\,\FF^{k-1}\right] \\
&= \|z^{k-1}-z^\star\|_2^2 + \langle z^{k-1} + y^k - 2 z^\star , \bP(y^k - z^{k-1} ) \rangle,
\end{align*}
where we used the fact that $z^{k-1}$ and $y^k$ are $\FF^{k-1}$-measurable and that $P_{\Sel^k}$ is a projection matrix so $P_{\Sel^k} = P_{\Sel^k}^\top = P_{\Sel^k}^2 $.

Then, using the fact $y^\star = z^\star$, the scalar product above can be simplified as follows
\begin{align*}
\nonumber\langle z^{k-1} &+ y^k - 2z^\star , \bP(y^k - z^{k-1}) \rangle = \langle z^{k-1} + y^k - z^\star - y^\star, \bP(y^k - z^{k-1} + y^\star -  z^\star ) \rangle \\
 \nonumber
 &=  - \langle z^{k-1}  - z^\star, \bP(z^{k-1} - z^\star) \rangle + \langle z^{k-1} - z^\star, \bP(y^k - y^\star) \rangle\\
 \nonumber
&\hspace*{0.2cm} + \langle y^k  - y^\star, \bP(y^k - y^\star) \rangle - \langle y^k - y^\star, \bP(z^{k-1} - z^\star) \rangle\\
&=  \langle y^k  - y^\star, \bP(y^k - y^\star) \rangle - \langle z^{k-1}  - z^\star, \bP(z^{k-1} - z^\star) \rangle
\end{align*}
where we used in the last equality that $\bP$ is symmetric.
 \endproof

\begin{lemma}[Contraction property in $\bP$-weighted norm]\label{lm:bub}
From the minimizer  $x^\star$  of \eqref{eq:main_problem}, define the fixed points $z^\star = y^\star = \bQ\left(x^\star - \gamma\nabla f\left(x^\star\right)\right)$ of the sequences $(y^k)$ and $(z^k)$. 
If Assumptions \ref{hyp:f} and \ref{hyp:main} hold, then
    \begin{equation*}
\|y^{k}-y^\star\|_{\bP}^2 -  \|z^{k-1}-z^\star\|_{\bP}^2 \leq - \lambda_{\min}(\bP) \frac{2\gamma\mu L}{\mu + L}\|z^{k-1} - z^\star\|_2^2.
\end{equation*}
\end{lemma}

\proof{Proof.}
First, using the definition of $y^k$ and $y^\star$,
\begin{align*}
  \|y^{k}-y^\star\|_{\bP}^2
   &=
   \langle \bQ(x^k - \gamma\nabla f(x^k) - x^\star + \gamma\nabla f(x^\star)), 
   \bP\bQ(x^k - \gamma\nabla f(x^k) - x^\star + \gamma\nabla f(x^\star)) \rangle\\
   &=
   \langle x^k - \gamma\nabla f(x^k) - x^\star + \gamma\nabla f(x^\star), \bQ^\top\bP\bQ (x^k - \gamma\nabla f(x^k) - x^\star + \gamma\nabla f(x^\star)) \rangle\\
   &= \left\| x^k - \gamma\nabla f(x^k) - ( x^\star - \gamma\nabla f(x^\star) ) \right\|_2^2.
\end{align*}
Using the standard stepsize range $\gamma\in(0,2/(\mu+L)]$, one has (see e.g. \cite[Lemma~3.11]{bubeck2015convex}) 
$$
 \|y^{k}-y^\star\|_{\bP}^2 = \left\| x^k - \gamma\nabla f(x^k) - ( x^\star - \gamma\nabla f(x^\star) ) \right\|_2^2 \leq \left(1 - \frac{2\gamma\mu L}{\mu + L}\right)\|x^k - x^\star\|_2^2.
$$
Using the non-expansiveness of the proximity operator of convex l.s.c. function $g$ \cite[Prop. ~12.27]{bauschke2011convex} along with the fact that as $x^\star$ is a minimizer of \eqref{eq:main_problem} so $x^\star = \prox_{\gamma g}(x^\star - \gamma \nabla f(x^\star)) = \prox_{\gamma g}(\bQ^{-1}z^\star)$ \cite[Th. ~26.2]{bauschke2011convex}, we get
\begin{align*}
\|x^k - x^\star\|_2^2 &= \|\prox_{\gamma g}(\bQ^{-1}(z^{k-1})) - \prox_{\gamma g}(\bQ^{-1}(z^\star))\|_2^2\\
&\leq
\|\bQ^{-1}(z^{k-1} - z^\star)\|_2^2 = \langle \bQ^{-1}(z^{k-1} - z^\star),\bQ^{-1}(z^{k-1} - z^\star)\rangle \\
&= \langle z^{k-1} - z^\star,\bP(z^{k-1} - z^\star)\rangle =  \|z^{k-1}-z^\star\|_{\bP}^2
\end{align*}
where we used that $\bQ^{-\top}\bQ^{-1} = \bQ^{-2} = \bP$. Combining the previous equations, we get 
\begin{align*}
    \|y^{k}-y^\star\|_{\bP}^2 -  \|z^{k-1}-z^\star\|_{\bP}^2 \leq- \frac{2\gamma\mu L}{\mu + L} \|z^{k-1} - z^\star\|_{\bP}^2.
\end{align*}
Finally, the fact that $\|x\|_{\bP}^2 \geq \lambda_{\min}(\bP)\|x\|_2^2$ for positive definite matrix $\bP$ enables to get the claimed result.
 \endproof

Relying on these two lemmas, we are now able to prove Theorem~\ref{th:conv_nondis}.
by showing that the distance of $z^k$ towards the minimizers is a contracting super-martingale. 

\proof{Proof.}[Proof of Theorem~\ref{th:conv_nondis}.]
Combining Lemmas~\ref{lm:removing_exp} and \ref{lm:bub}, we get
\begin{align*}
    \EE\left[\|z^{k} - z^\star\|_2^2\,|\,\FF^{k-1}\right] \leq \left(1 - \lambda_{\min}(\bP)\frac{2\gamma\mu L}{\mu + L}\right)\|z^{k-1}-z^\star\|_2^2
\end{align*}
and thus by taking the full expectation and using nested filtrations $(\FF^k)$, we obtain
\begin{equation*}
    \EE\left[\|z^{k} - z^\star\|_2^2\right] \leq \left(1 - \lambda_{\min}(\bP)\frac{2\gamma\mu L}{\mu + L}\right)^k\|z^0-z^\star\|_2^2 = \left(1 - \lambda_{\min}(\bP)\frac{2\gamma\mu L}{\mu + L}\right)^k\|z^0-\bQ(x^\star - \gamma\nabla f(x^\star))\|_2^2.
\end{equation*}
Using the same arguments as in the proof of Lemma~\ref{lm:bub}, one has
\begin{align*}
  \|x^{k+1}-x^\star\|_2^2 \leq  \|z^k-z^\star\|_{\bP}^2 \leq \lambda_{\max}(\bP)\|z^k-z^\star\|_2^2
\end{align*}
which enables to conclude 
\begin{equation*}
    \EE\left[\|x^{k+1}-x^\star\|_2^2\right] \leq \left(1 - \lambda_{\min}(\bP)\frac{2\gamma\mu L}{\mu + L}\right)^k \lambda_{\max}(\bP) \|z^0-\bQ(x^\star - \gamma\nabla f(x^\star))\|_2^2.
\end{equation*}
Finally, this linear convergences implies the almost sure convergence of $(x^k)$ to $x^\star$ as 
\begin{equation*}
    \EE\left[ \sum_{k=1}^{+\infty}  \|x^{k+1} - x^\star \|^2 \right]\leq C \sum_{k=1}^{+\infty}\left(1-\lambda_{\min}(\bP)\frac{2\gamma\mu L}{\mu + L}\right)^k < +\infty
\end{equation*}
implies that $ \sum_{k=1}^{+\infty}  \|x^{k+1} - x^\star \|^2$ is finite with probability one. Thus we get
\[
1 = \mathbb{P}\left[ \sum_{k=1}^{+\infty}  \|x^{k+1} - x^\star \|^2 < +\infty \right] \leq \mathbb{P}\left[ \|x^k - x^\star \|^2 \to 0 \right]
\]
which in turn implies that $(x^k)$ converges almost surely to $x^\star$.
 \endproof

\subsection{Examples and connections with the existing work}\label{sec:comparison}
 
In this section, we derive specific cases and discuss the relation between our algorithm and the related literature.

\subsubsection{Projections onto coordinates} \label{sec:coordproj}

A simple instantiation of our setting can be obtained by considering projections onto uniformly chosen coordinates (Example\;\ref{ex:P}); with the family 
$$ \C = \{\C_1,..,\C_n\} \quad\text{ with } \C_i = \{x\in\mathbb{R}^n  : x_j = 0 ~~ \forall j \neq i \} $$ 
and the selection $\Sel$ consisting of taking $\C_i$ according to the output of a Bernoulli experiment of parameter $p_i$. Then, the matrices $\bP = \diag( [p_1,..,p_n])$, $P_{\Sel^k}$ and $\bQ$ commute, and, by a change of variables $\tilde{z}^k = \mathsf{Q}^{-1} z^k$ and $\tilde{y}^k  = \mathsf{Q}^{-1} y^k$, Algorithm \ref{alg:strata_nondis} boils down to
\[
\tilde{y}^k = x^k - \gamma\nabla f\left(x^k\right)\qquad
\tilde{z}^{k} =  P_{\Sel^k}\left(\tilde{y}^k\right) + (I- P_{\Sel^k})\left(\tilde{z}^{k-1}\right),\qquad
x^{k+1}  = \prox_{\gamma g} \left(\tilde{z}^{k}\right)
 \]
i.e. no change of basis is needed anymore, even if $g$ is non-separable. Furthermore, the convergence rates simplifies to $(1- 2 \min_i p_i  \gamma \mu L/(\mu + L))$ which translates to $(1- 4 \min_i p_i\mu L/(\mu + L)^2)$ for the optimal $\gamma = 2/(\mu+L)$.

In the special case where $g$ is separable (i.e. $g(x) = \sum_{i=1}^n g_i(x_i)$), we can further simplify the iteration. In this case, projection and proximal steps commute, so that the iteration can be written
\begin{align*}
 x^{k+1} &=  P_{\Sel^k}\prox_{\gamma g} \left( x^k - \gamma\nabla f(x^k)\right)  + (I -  P_{\Sel^k}) x^k\\
\text{i.e. } x_i^{k+1} &= \left\{  \begin{array}{ll}
   \prox_{\gamma g_i} \left( x_i^k - \gamma\nabla_i f(x^k)\right) = \displaystyle\arg\min_w g_i(w) + \langle w , \nabla_i f(x^k) \rangle  + \frac{1}{2\gamma}\|w - x_i^k\|_2^2 & \text{ if } i \in \Sel^k  \\
   x_i^k  & \text{ elsewhere}
\end{array}  \right. 
\end{align*}
which boils down to the usual (proximal) coordinate descent algorithm, that recently knew a rebirth in the context of huge-scale optimization, see\;\cite{tseng2001convergence}, \cite{nesterov2012efficiency}, \cite{richtarik2014iteration} or \cite{wright2015coordinate}.
In this special case, the theoretical convergence rate of \algo is close to the existing rates in the literature. 
For clarity, 
we compare with the uniform randomized coordinate descent of \cite{richtarik2014iteration} 
(more precisely Th.\;6 with $L_i =L$, $B_i=1$, $\mu L\leq 2$) which can be written as $\left(1 - \mu L/4n\right)$ in $\ell_2$-norm.
The rate of \algo in the same uniform setting (Example~\ref{ex:P} with $p_i = p=1/n$) is $\left(1 - \frac{4\mu L}{n(\mu+ L)^2}\right)$ with the optimal step-size.

\subsubsection{Projections onto vectors of fixed variations}
\label{sec:var}

The vast majority of randomized subspace methods consider the coordinate-wise projections treated in \ref{sec:coordproj}. This success is notably due to the fact that most problems onto which they are applied have naturally a coordinate-wise structure; for instance, due to the structure of $g$ ($\ell_1$-norm, group lasso, etc). However, many problems in signal processing and machine learning feature a very different structure. A typical example is when $g$ is the $1$D-Total Variation 
\begin{equation}\label{eq:TV}
    g(x) = \sum_{i=2}^n |x_{i}-x_{i-1}|
\end{equation}
featured for instance in the fused lasso problem \cite{tibshirani2005sparsity}. In order to project onto subspaces of vectors of fixed variation (i.e. vectors for which $x_j = x_{j+1}$ except for a prescribed set of indices), one can define the following covering family 
$$ \C = \{\C_1,..,\C_{n-1}\} \quad\text{ with } \C_i = \left\{x\in\mathbb{R}^n  : x_{j} = x_{j+1} \text{ for all } j\in \{1,..,n-1\}\setminus\{i\}  \right\} $$
and an admissible selection $\Sel$ consisting in selecting uniformly $s$ elements in $\C$. 
Then, if $\Sel$ selects $\C_{n_1},...,\C_{n_s}$, the update will live in the sum of these subspaces, i.e. the subspace of the vectors having jumps at coordinates $n_1,n_2,..,n_s$. Thus, the associated projection in the algorithm writes\\[2ex]
{\footnotesize 
\begin{equation}\label{eq:proj_tv}
P_{\Sel} =
\begin{matrix}
\begin{pmatrix}
\overmat{$n_1$}{\frac{1}{n_1}& \hdots & \frac{1}{n_1}}& 0 &\hdots & \overmat{$n-n_s$}{\hdots &\hdots & \hspace*{4pt}0\hspace*{4pt}} \\
\vdots& \ddots & \vdots & \vdots & \ddots & \ddots &\ddots & \vdots\\
\frac{1}{n_1}& \hdots & \frac{1}{n_1}& 0 &\ddots & \ddots &\ddots & \vdots\\
0 & \hdots & 0  &\ddots & \ddots & \ddots & \ddots & \vdots\\
\vdots & \ddots & \ddots & \ddots  &\ddots & 0  & \hdots & 0\\
\vdots &\ddots & \ddots &\ddots & 0 & \frac{1}{n-n_s}& \hdots & \frac{1}{n-n_s}\\
\vdots & \ddots & \ddots &\ddots & \vdots & \vdots& \ddots & \vdots\\
0 &\hdots & \hdots &\hdots & 0 & \frac{1}{n-n_s}& \hdots & \frac{1}{n-n_s}
\end{pmatrix}
\hspace*{-5pt}
\begin{aligned}
  &\left.\begin{matrix}
  \\[36pt]
  \end{matrix} \right\} %
  n_1\\[40pt]
  &\left.\begin{matrix}
  \\[36pt]
  \end{matrix} \right\} %
  n-n_s
\end{aligned}
\end{matrix}
\end{equation}}
Note also that $P_{\Sel}x$ has the same value for coordinates $[n_i,n_{i+1})$, equal to the average of these values.

As mentioned above, the similarity between the structure of the optimization problem and the one of the subspace descent is fundamental for performance in practice. In Section\;\ref{sec:adapt}, we exploit the identification properties of the proximity operator in order to automatically adapt the subspace selection, which leads to a tremendous gain in performance.

\subsubsection{Comparison with sketching}

In sharp contrast with the existing literature, our subspace descent algorithm handles non-separable regularizers $g$. A notable exception is the algorithm called \sega \cite{hanzely2018sega}, a random sketch-and-project proximal algorithm, that can also deal with non-separable regularizers. While the algorithm shares similar components with ours, the main differences between the two algorithms are 
\begin{itemize}
    \item {biasedness of the gradient}: \sega deals with unbiased gradients while they are biased for \algo;
\item {projection type}: \sega projects the gradient while we project after a gradient step (option (b) vs.\,option (c) in the discussion starting Section\;\ref{sec:algo}).
\end{itemize} 
These differences are fundamental and create a large gap in terms of target, analysis and performance between the two algorithms. The practical comparison is illustrated in Section\;\ref{sec:num:sega}.

\section{Adaptive subspace descent}

\label{sec:ada_sub_des}

This section presents an extension of our randomized subspace descent algorithm where the projections are iterate-dependent. 
Our aim is to automatically adapt to the structure identified by the iterates along the run of the algorithm.

The methods proposed here are, up to our knowledge, the first ones where the iterate structure enforced by a nonsmooth regularizer is used to adapt the selection probabilities in a randomized first-order method. As discussed in the introduction, even for the special case of coordinate descent, our approach is different from existing techniques that use fixed arbitrary probabilities \cite{richtarik2014iteration,necoara2014random}, greedy selection \cite{dhillon2011nearest,nutini2015coordinate,nutini2017let}, or adaptive selection 
based on the coordinate-wise Lipschitz constant and coordinates \cite{perekrestenko2017faster,namkoong2017adaptive,stich2017safe}. 

We present our adaptive subspace descent algorithm in two steps. First, we introduce in Section~\ref{sec:ada_algo} a generic algorithm with varying selections and establish its convergence. 
Second, in Section~\ref{sec:identif}, we provide a simple general identification result. We then combine these two results to provide an efficient adaptive method in Section~\ref{sec:adapt}. 

\subsection{Random subspace descent with time-varying selection}
\label{sec:ada_algo}

For any randomized algorithm, using iterate-dependent sampling would  automatically break down the usual i.i.d. assumption.
In our case, adapting to the current iterate structure means that the associated random variable depends on the past. We thus need further analysis and notation.

In the following, we use the subscript $\ell$ to denote the $\ell$-th change in the selection. We denote by $\Up$ the set of time indices at which an adaptation is made, themselves denoted by $k_\ell  = \min \{ k> k_{\ell-1} : k\in\Up\}$. 

In practice, at each time $k$, there are two decisions to make (see Section\;\ref{sec:adapt}): (i) \emph{if} an adaptation should be performed; and (ii) \emph{how} to update the selection. Thus, we replace the i.i.d. assumption of Assumption\;\ref{hyp:main} with the following one.

\begin{assumption}[On the randomness of the adaptive algorithm]\label{hyp:main_identif}
For all $k>0$, $\Sel^k$ is $\FF^k$-measurable and admissible. Furthermore, if $k\notin\Up$,  $(\Sel^k)$ is independent and identically distributed on $[{k}_{\ell},k]$. The decision to adapt or not at time $k$ is $\FF^k$-measurable, i.e. $(k_\ell)_\ell$ is a sequence of $\FF^k$-stopping times. 
\end{assumption}

Under this assumption, we can prove the convergence of the varying-selection random subspace descent, Algorithm \ref{alg:ada_strata_nondis}. A generic result  is given in Theorem\;\ref{th:conv_nondis_arbitrary} and a simple specification in the following example. The rationale of the proof is that the stability of the algorithm is maintained when adaptation is performed sparingly.

\begin{algorithm}[H] 
\caption{Adaptive Randomized Proximal Subspace Descent - \adaalgo}
\label{alg:ada_strata_nondis}
\begin{algorithmic}[1] 
    \STATE Initialize $z^0$, $x^1 = \prox_{\gamma g}(\bQ_0^{-1}(z^0))$, $\ell=0$, $\Up=\{0\}$.
    \FOR{$k=1,\ldots$}
            \STATE $y^k = \bQ_{\ell}\left(x^k - \gamma\nabla f\left(x^k\right)\right)$
            \STATE $z^{k} = P_{\Sel^k} \left(y^k\right) + (I- P_{\Sel^k} ) \left(z^{k-1}\right)$
            \STATE$x^{k+1} = \prox_{\gamma g} \left(\bQ_\ell^{-1}\left(z^{k}\right)\right)$
            \IF{an adaptation is decided }
            \STATE $\Up \leftarrow \Up \cup \{k+1\}$, $\ell\leftarrow \ell +1$
            \STATE Generate a new admissible selection
            \STATE Compute $\bQ_\ell = \bP_\ell^{-\frac12}$ and $\bQ_\ell^{-1}$
            \STATE Rescale $z^k \leftarrow \bQ_\ell   \bQ_{\ell-1}^{-1} z^k$ \label{line:rescale}
            \ENDIF
    \ENDFOR
\end{algorithmic}
\end{algorithm}

\begin{theorem}[\adaalgo~convergence]\label{th:conv_nondis_arbitrary}
Let Assumptions \ref{hyp:f} and  \ref{hyp:main_identif} hold.  For any $\gamma\in(0,2/(\mu+L)]$, let the user choose its adaptation strategy so that:
\begin{itemize}
    \item  the \emph{adaptation cost} is upper bounded by a deterministic sequence:  $ \|  \bQ_\ell   \bQ_{\ell-1}^{-1} \|_2^2 \leq \mathbf{a}_\ell $;
    \item the \emph{inter-adaptation time} is lower bounded by a deterministic sequence: $k_{\ell}-k_{\ell-1}\geq \mathbf{c}_\ell$;
    \item the \emph{selection uniformity} is lower bounded by a deterministic sequence: $\lambda_{ \min}(\bP_{\ell-1}) \geq \lambda_{\ell-1} $;
\end{itemize}
then, from the \emph{previous instantaneous rate} $1-\alpha_{\ell-1}  := 1 -  2\gamma \mu L \lambda_{\ell-1}/(\mu + L)  $, the \emph{corrected rate} for cycle $\ell$ writes  
\begin{equation}\label{eq:corr-rate}
(1-\beta_\ell) := (1-\alpha_{\ell-1})\mathbf{a}_\ell^{1/\mathbf{c}_\ell}. 
\end{equation}
Then, we have for any $k\in [k_\ell,k_{\ell+1})$
\begin{align*}
    \EE\left[\|x^{k+1}-x^\star\|_2^2\right] &\leq (1-\alpha_\ell)^{k-k_\ell} \prod_{m=1}^\ell (1-\beta_m)^{\mathbf{c}_m}  \|z^{0}-\bQ_0\left(x^\star - \gamma\nabla f\left(x^\star\right)\right)\|_2^2.
\end{align*}
\end{theorem}

This theorem means that by balancing the magnitude of the adaptation (i.e.\;$\mathbf{a}_m$) with the time before adaptation (i.e.\;$\mathbf{c}_m$) from the knowledge of the current rate $(1-\alpha_{m-1})$, one can retrieve the exponential convergence with a controlled degraded rate $(1-\beta_m)$. 
This result is quite generic, but it can be easily adapted to specific situations.
For instance, we provide a simple example with a global rate on the iterates in the forthcoming Example\;\ref{ex:adapt_fixed}.

\begin{figure}[!ht]
    \centering

\resizebox {0.9\columnwidth} {!} {
     \begin{tikzpicture}
    \draw[thick, ->] (0,0) -- (10,0) node [below,scale=0.7] {iterations};
    \foreach \x in {1,...,4}
    \draw (2*\x, 0.1) -- node[pos=0.5] (point\x) {} (2*\x, -0.1);

          \draw[thick,dashed,white] (2.5,0) -- (3.5,0);
          \draw[thick,dashed,white] (4.5,0) -- (5.5,0);
    
    \node[scale=0.7] at (2,0.25) {\color{blue!70!black} $k_\ell$};

    \draw[thick,blue!70!black, ->] (2,0.7) -- (2,0.4) ;
    \node[scale = 0.7] at (1.9,0.9) {\color{blue!70!black} adaptation};


    
        \draw[thick,green!30!black, ->] (4,0.7) -- (4,0.4) ;
    \node[scale = 0.7] at (4,0.25) {\color{green!30!black} $k_\ell + \mathbf{c}_{\ell+1}$};
        \node[scale = 0.7] at (4.1,0.9) {\color{green!30!black} new adaptation possible};

    \node[scale = 0.5, gray] at (5.2,-0.5) { $z^{k-1},x^{k},y^{k}$};
    
    \node[scale = 0.7] at (6,0.25) { $k$};
    \node[scale = 0.7] at (6,-0.25) { $\Sel^k$};

    \draw[gray] (5.8,-0.15) -- (5.8,-0.7) -- (3,-0.7);
    \draw[dashed,gray] (2.5,-0.7) -- (3,-0.7);
    \node[scale = 0.7,gray] at (4,-0.6) { $\FF^{k-1}$};

    \node[scale = 0.5] at (7,-0.5) { $z^k\to x^{k+1} \to y^{k+1}$};
    \node[scale = 0.5] at (7,-0.7) { $\{ k_{\ell+1} = {k+1}\}$};
    
    \draw (7.7,-0.15) -- (7.7,-1.0) -- (3,-1.0);
    \draw[dashed] (2.5,-1.0) -- (3,-1.0);
    \node[scale = 0.7] at (4,-0.9) { $\FF^{k}$};

    \node[scale = 0.7] at (8,0.25) { $k+1$};
    \node[scale = 0.7,red] at (8,-0.25) { $\Sel^{k+1}$};

    

  \end{tikzpicture}}

      \caption{Summary of notations about iteration, adaptation and filtration. The filtration $\mathcal{F}^{k-1}$ is the sigma-algebra generated by $\{\Sel^\ell\}_{\ell\leq k-1}$ encompassing the knowledge of all variables up to $y^k$ (but not $z^k$).}
    \label{fig:proof}
  
  \end{figure}
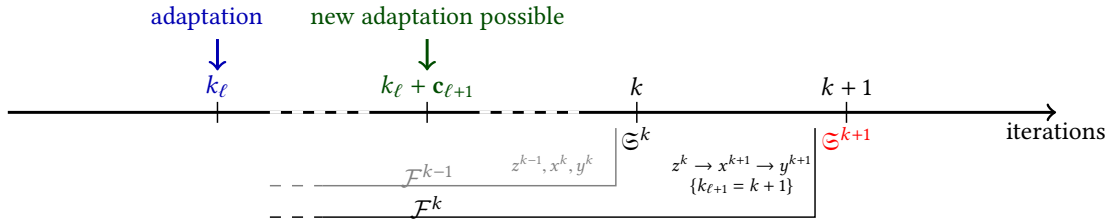

For now, let us turn to the proof of the theorem. To ease its reading, the main notations and measurability relations are depicted in Figure\;\ref{fig:proof}.

\proof{Proof.}
We start by noticing that, for a solution $x^\star$ of \eqref{eq:main_problem}, the proof of Theorem~\ref{th:conv_nondis} introduces the companion variable $z^\star =  \bQ\left(x^\star - \gamma\nabla f\left(x^\star\right)\right)$ which directly depends on $\bQ$, preventing us from a straightforward use of the results of Section\;\ref{sec:conv}. However, defining $z^{\star}_{\ell} =  \bQ_\ell\left(x^\star - \gamma\nabla f\left(x^\star\right)\right)$, Lemmas~\ref{lm:removing_exp} and \ref{lm:bub} can be directly extended and combined to show for any $k\in [k_\ell,k_{\ell+1})$ 
\begin{equation}\label{eq:iterate_flexible_lambda}
    \EE\left[\|z^{k} - z^{\star}_{\ell}\|_2^2\,|\,\FF^{k-1}\right] \leq \underbrace{  \left(1 -  \frac{2\gamma \mu L \lambda_{ \min}(\bP_\ell)}{\mu + L} \right) }_{\leq 1- \alpha_\ell} \|z^{k-1}-z^{\star}_{\ell}\|_2^2 .
\end{equation} 
Since the distribution of the selection has not changed since ${k}_{\ell}$, iterating \eqref{eq:iterate_flexible_lambda} leads to 
\begin{align}\label{eq:iterate_flexible_lambda2}
    \EE\left[\|z^{k} - z^{\star}_{\ell}\|_2^2\,|\,\FF^{k_{\ell}-1}\right] &\leq (1- \alpha_\ell)^{k-k_\ell} \|z^{k_\ell-1}-z^{\star}_{\ell}\|_2^2 .
\end{align}
We focus now on the term $\|z^{k_\ell-1}-z^{\star}_{\ell}\|_2^2$ corresponding to what happens at the last adaptation step. From the definition of variables in the algorithm and using the deterministic bound on $\|  \bQ_\ell   \bQ_{\ell-1}^{-1}\|$, we write
\begin{align}
  \nonumber  \EE\left[\|z^{k_\ell-1}-z^{\star}_{\ell}\|_2^2\,|\,\FF^{k_{\ell}-2}\right] &\leq   \EE\left[ \|  \bQ_\ell   \bQ_{\ell-1}^{-1}(z^{k_{\ell}-2} + P_{k_{\ell}-1}(y^{k_{\ell}-1}-z^{k_{\ell}-2}) - \bQ_\ell   \bQ_{\ell-1}^{-1}z^{\star}_{\ell-1}\|_2^2\,|\,\FF^{k_{\ell}-2}\right] \\[1.5ex]
        &\leq  \EE\left[ \|  \bQ_\ell   \bQ_{\ell-1}^{-1} \|_2^2   \|  z^{k_{\ell}-2} + P_{k_{\ell}-1}(y^{k_{\ell}-1}-z^{k_{\ell}-2}) - z^{\star}_{\ell-1}\|_2^2\,|\,\FF^{k_{\ell}-2}\right] \label{eq:Jadd} \\[1.5ex]
   \nonumber &\leq   \mathbf{a}_\ell  (1-\alpha_{\ell-1}) \|z^{k_\ell-2}-z^{\star}_{\ell-1}\|_2^2. \label{eq:adapt}
\end{align} 
Repeating this inequality backward to the previous adaptation step $z^{k_{\ell-1}}$, we get 
\begin{align}
  \nonumber  \EE\left[\|z^{k_\ell-1}-z^{\star}_{\ell}\|_2^2\,|\,\FF^{k_{\ell-1}}\right] &\leq  \mathbf{a}_\ell  (1-\alpha_{\ell-1})^{k_{\ell}-k_{\ell-1}} \|z^{k_{\ell-1}}-z^{\star}_{\ell-1}\|_2^2 \\
   &\leq  \mathbf{a}_\ell  (1-\alpha_{\ell-1})^{\mathbf{c}_\ell} \|z^{k_{\ell-1}}-z^{\star}_{\ell-1}\|_2^2,
\end{align} 
using the assumption of bounded inter-adaptation times.
Combining this inequality and \eqref{eq:iterate_flexible_lambda2}, we obtain that for any $k\in[k_\ell,k_{\ell+1})$, 
\begin{align*}
    \EE\left[\|z^{k}-z^{\star}_{\ell}\|_2^2\right] &\leq (1-\alpha_\ell)^{k-k_\ell}
    \prod_{m=1}^\ell  \mathbf{a}_m (1-\alpha_{m-1})^{\mathbf{c}_m} \|z^{0}-z^{\star}_{0}\|_2^2.
    \end{align*}
Using now \eqref{eq:corr-rate}, we get    
\[
\EE\left[\|z^{k}-z^{\star}_{\ell}\|_2^2\right] ~\leq~ (1-\alpha_\ell)^{k-k_\ell}
    \prod_{m=1}^\ell  (1-\beta_{m})^{\mathbf{c}_m} \|z^{0}-z^{\star}_{0}\|_2^2
   \]
Finally, the non-expansiveness of the prox-operator propagates this inequality to $x_k$, since we have 
\begin{align*}
\|x^k - x^\star\|_2^2 &= \|\prox_{\gamma g}(\bQ_\ell^{-1}(z^{k-1})) - \prox_{\gamma g}(\bQ_\ell^{-1}(z_\ell^\star))\|_2^2\\
&\leq
\|\bQ_\ell^{-1}(z^{k-1} - z_\ell^\star)\|_2^2 \leq \lambda_{\max}(\bQ_\ell^{-1})^2  \|z^{k-1}-z^\star_\ell\|_2^2 =  
\lambda_{\max}(\bP_\ell)  \|z^{k-1}-z^\star_\ell\|_2^2 \leq \|z^{k-1}-z^\star_\ell\|_2^2.
\end{align*}
This concludes the proof.
 \endproof

\begin{example}[Explicit convergence rate] \label{ex:adapt_fixed}
Let us specify Theorem~\ref{th:conv_nondis_arbitrary} with the following simple adaptation strategy. We take a fixed upper bound on the adaptation cost and a fixed lower bound on uniformity:
\begin{equation}\label{eq:explicit}
\|  \bQ_\ell   \bQ_{\ell-1}^{-1} \|_2^2 \leq \mathbf{a} \qquad
\lambda_{ \min}(\bP_{\ell}) \geq \lambda.
\end{equation}
Then from the rate $1-\alpha = 1- 2\gamma\mu L\lambda/(\mu+L)$, we can perform an adaptation every 
\begin{align}
    \label{eq:min_adapt}
\mathbf{c} = \lceil \log(\mathbf{a})/\log\big((2-\alpha)/(2-2\alpha)\big)\rceil
\end{align}
 iterations, so that $\mathbf{a}(1-\alpha)^\mathbf{c} = (1-\alpha/2)^\mathbf{c}$ and $k_\ell = \ell \mathbf{c}$. A direct application of Theorem \eqref{th:conv_nondis_arbitrary} gives that, for any $k$,
 \begin{align*}
    \EE\left[\|x^{k+1}-x^{\star}_{\ell}\|_2^2\right] &\leq  \left(1-\frac{\gamma\mu L\lambda}{\mu+L}\right)^{k}  C
\end{align*}
where $C = \|z^0 - \bQ_0(x^\star - \gamma\nabla f(x^\star))\|_2^2$. That is the same convergence mode as in the non-adaptive case (Theorem\;\ref{th:conv_nondis}) with a modified rate. Note the modified rate provided here (of the form $(1-\alpha/2)$ to be compared with the $1-\alpha$ of Theorem\;\ref{th:conv_nondis}) was chosen for clarity; any rate strictly slower than $1-\alpha$ can bring the same result by adapting $\mathbf{c}$ accordingly. 
\end{example}

\begin{remark}[On the adaptation frequency]
Theorem\;\ref{th:conv_nondis_arbitrary} and Example\;\ref{ex:adapt_fixed} tell us that we have to respect a prescribed number of iterations between two adaptation steps. We emphasize here that if this inter-adaptation time is violated, the resulting algorithm may be highly unstable. We illustrate this phenomenon on a TV-regularized least squares problem: we compare two versions of \adaalgo with the same adaptation strategy verifying 
\eqref{eq:explicit} but with two different adaptation frequencies
\begin{itemize}
    \item at every iteration (i.e. taking $\mathbf{c}_\ell = 1$)
    \item following theory (i.e. taking $\mathbf{c}_\ell = \mathbf{c}$ as per Eq.~\eqref{eq:min_adapt})
\end{itemize}
On Figure~\ref{fig:stab}, we observe that adapting every iteration leads to a chaotic behavior. Second, even though the theoretical number of iterations in an adaptation cycle is often pessimistic (due to the rough bounding of the rate), the iterates produced with this choice quickly become stable (i.e. identification happens, which will be shown and exploited in the next section) and show a steady decrease in suboptimality.   

\begin{figure}[H]
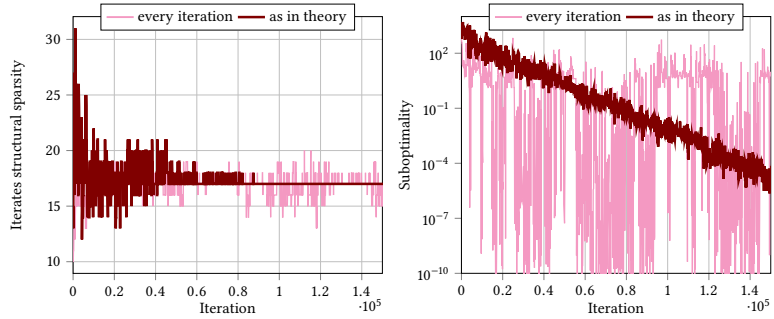

\begin{center}
 \scalebox{.6}{\input{Article/fused_lasso_sparsity.tex}}
 \scalebox{.6}{\input{Article/fused_lasso_suboptimality_vs_iteration.tex}}
\end{center}
\caption{Comparisons between theoretical and harsh updating time for \adaalgo.}
\label{fig:stab}
\end{figure}

\end{remark}

A drawback of Theorem~\ref{th:conv_nondis_arbitrary} is that the adaptation cost, inter-adaptation time, and selection uniformity have to be bounded by deterministic sequences. This can be restrictive if we do not have prior knowledge on the problem or if the adaptation cost varies a lot. 
This drawback can be circumvented to the price of loosing the rate \emph{per iteration} to the rate \emph{per adaptation}, as formalized in the following result.

\begin{theorem}[\adaalgo convergence: practical version]\label{th:aggressive}
    Let Assumptions \ref{hyp:f} and  \ref{hyp:main_identif} hold. Take $\gamma\in(0,2/(\mu+L)]$, choose $\lambda>0$, and set $\beta = \gamma\mu L\lambda/(\mu+L)$. Consider the following adaptation strategy:
    \begin{itemize}
        \item[1)] From the observation of $x^{k_{\ell-1}}$\!,\;choose a new sampling with $\bP_{\ell}$ and $\bQ_\ell$,\,such\;that\;$\lambda_{ \min}(\bP_{\ell}) \geq \lambda$;
        \item[2)] Compute $\mathbf{c}_\ell$ so that $\|  \bQ_\ell   \bQ_{\ell-1}^{-1} \|_2^2 (1-\alpha_{\ell-1})^{\mathbf{c}_\ell} \leq 1 - \beta$ where $\alpha_{\ell-1} = 2\gamma \mu L \lambda_{ \min}(\bP_{\ell-1})/(\mu + L)$;
        \item[3)] Apply the new sampling 
        after $\mathbf{c}_\ell$ iterations ($k_{\ell} = k_{\ell-1}+\mathbf{c}_\ell$).
    \end{itemize}
    Then, we have for any $k\in [k_\ell,k_{\ell+1})$
    \begin{align*}
        \EE\left[\|x^{k+1}-x^\star\|_2^2\right] &\leq (1-\alpha_\ell)^{k-k_\ell}  \left(1-\beta\right)^{\ell}  \|z^{0}-\bQ_0\left(x^\star - \gamma\nabla f\left(x^\star\right)\right)\|_2^2. 
    \end{align*}
    \end{theorem}

\proof{Proof.}
The proof follows the same pattern as the one of Theorem~\ref{th:conv_nondis_arbitrary}.
The only difference is that the three control sequences (adaptation cost, inter-adaptation time, and selection uniformity) are now random sequences
since they depend on the iterates of the (random) algorithm. This technical point requires a special attention. 
In \eqref{eq:Jadd}, the adaptation introduces a cost by a factor $\|  \bQ_\ell   \bQ_{\ell-1}^{-1} \|_2^2 $, which is not deterministically upper-bounded  anymore. However it is $\FF^{k_{\ell-1}}$-measurable by construction of $\bQ_\ell$, so we can write
\begin{align*}
  \nonumber & \EE\left[\|z^{k_\ell-1}-z^{\star}_{\ell}\|_2^2\,|\,\FF^{k_{\ell-1}}\right] \\[1.5ex]
  \nonumber & = \EE\left[\EE\left[\|z^{k_\ell-1}-z^{\star}_{\ell}\|_2^2\,|\,\FF^{k_{\ell}-2}\right]\,|\,\FF^{k_{\ell-1}}\right]\\[1.5ex]
  &\leq   \EE\left[\EE\left[ \|  \bQ_\ell   \bQ_{\ell-1}^{-1}(z^{k_{\ell}-2} + P_{k_{\ell}-1}(y^{k_{\ell}-1}-z^{k_{\ell}-2}) - \bQ_\ell   \bQ_{\ell-1}^{-1}z^{\star}_{\ell-1}\|_2^2\,|\,\FF^{k_{\ell}-2}\right]\,|\,\FF^{k_{\ell-1}}\right]\\[1.5ex]
   \nonumber     &\leq \EE\left[\|  \bQ_\ell   \bQ_{\ell-1}^{-1} \|_2^2  (1-\alpha_{\ell-1}) \|z^{k_\ell-2}-z^{\star}_{\ell-1}\|_2^2\,|\,\FF^{k_{\ell-1}}\right]\\[1.5ex]
      \nonumber     &= \| \bQ_\ell   \bQ_{\ell-1}^{-1} \|_2^2  (1-\alpha_{\ell-1}) \EE\left[ \|z^{k_\ell-2}-z^{\star}_{\ell-1}\|_2^2\,|\,\FF^{k_{\ell-1}}\right].
   \end{align*} 
   Using Eq.~\eqref{eq:iterate_flexible_lambda}, this inequality yields  
  \begin{align*}
   \EE\left[\|z^{k_\ell-1}-z^{\star}_{\ell}\|_2^2\,|\,\FF^{k_{\ell-1}}\right]
   &\leq \|  \bQ_\ell   \bQ_{\ell-1}^{-1} \|_2^2  (1-\alpha_{\ell-1})^{k_\ell-k_{\ell-1}}  \EE\left[  \|z^{k_{\ell-1}-1}-z^{\star}_{\ell-1}\|_2^2\,|\,\FF^{k_{\ell-1}}\right]\\[1ex]
&\leq (1-\beta) \EE\left[  \|z^{k_{\ell-1}-1}-z^{\star}_{\ell-1}\|_2^2\,|\,\FF^{k_{\ell-1}}\right].
\end{align*}
where we used points $2)$ and $3)$ of the strategy to bound the first terms deterministically. Finally, we obtain
\begin{align*}
\EE\left[\|z^{k_\ell-1}-z^{\star}_{\ell}\|_2^2\right]
& = \EE\left[\EE\left[\|z^{k_\ell-1}-z^{\star}_{\ell}\|_2^2\,|\,\FF^{k_{\ell-1}}\right]\right]\\[1.5ex]
&\leq (1-\beta) \EE\left[\|z^{k_{\ell-1}-1}-z^{\star}_{\ell-1}\|_2^2\right]
\end{align*}
then the rest of the proof follows directly by induction.
 \endproof

\subsection{Identification of proximal algorithms}\label{sec:identif}

As discussed in the introduction, identification of some optimal structure has been extensively studied
in the context of constrained convex optimization (see e.g.\;\cite{wright1993identifiable}) and nonsmooth optimization (see e.g.\;\cite{hare2004identifying}).
In this section, we provide a general identification result for proximal algorithms useful for our developments, using the notion of sparsity vector.

\begin{definition}[Sparsity vector]\label{def:sparsity}
Let $\M= \{ \M_1,\ldots,\M_m\}$ be a family of subspaces of $\mathbb{R}^n$ with $m$ elements. We define the {sparsity vector} on $\M$ for point $x\in\mathbb{R}^n$ as the $\{0,1\}$-valued\footnote{For two vectors $a,b\in\{0,1\}^m$, we use the following notation and terminology: (1) if  $[a]_i \leq [b]_i$ for all $i=1,..,m$, we say that $b$ is greater than $a$, noted $a\leq b$; and (2) we define the union $c = a\cup b$ as $[c]_i = 1 $ if $[a]_i = 1$ or $[b]_i=1$ and $0$ elsewhere.} vector  $\S(x)\in \{0,1\}^m$ verifying
\begin{align}\label{eq:strata}
    \left[ \S(x) \right]_i = 0 \quad\text{ if } x \in \M_i \text{ and } 1 \text{ elsewhere}.
\end{align}
\end{definition}

An identification result is a theorem stating that the iterates of the considered algorithm eventually belong to some -- but not all -- subspaces in $\M$. We formulate such a result for almost surely converging proximal-based algorithms as follows. This very simple result is inspired from the extended identification result of \cite{fadili2018sensitivity} (but does not rely on strong primal-dual structures as presented in \cite{fadili2018sensitivity}).

\begin{theorem}[Enlarged identification]\label{th:identification}
Let $(u^k)$ be an $\mathbb{R}^n$-valued sequence converging almost surely to $u^\star$ and define sequence $(x^k)$ as $x^k = \prox_{\gamma g}(u^k)$ and  $x^\star =  \prox_{\gamma g}(u^\star)$. Then $(x^k)$ identifies some subspaces with probability one; more precisely for any $\varepsilon>0$, with probability one, after some finite time,
\begin{equation}\label{eq:identification_result}
	 \S(x^\star) ~\leq~ \S(x^k) ~\leq\!\!\bigcup_{u\in\mathcal{B}(u^\star,\varepsilon)} \!\S(\prox_{\gamma g}(u)).
\end{equation}
\end{theorem}

\proof{Proof.}
The proof is divided between the two inequalities. We start with the right inequality. As $u^k \to u^\star$ almost surely, for any $\varepsilon>0$, $u^k$ will belong to a ball centered around $u^\star$ of radius $\varepsilon$ in finite time with probability one. Then, trivially, it will belong to a subspace if all points in this ball belong to it, which corresponds to the second inequality.

Let us turn now to the proof of the left inequality.  Consider the sets to which $x^\star$ belongs i.e. $\M^\star = \{ \M_i\in\M : x^\star \in \M_i \}$; as $\M$ is a family of subspaces, there exists a ball of radius $\varepsilon'>0$ around $x^\star$ such that no point $x$ in it belong to {more} subspaces than $x^\star$ i.e. $x \notin \M \setminus \M^\star$. As $x^k \to x^\star $ almost surely, it will reach this ball in finite time with probability one and thus belong to fewer subspaces than $x^\star$.
 \endproof

This general theorem explains that iterates of any converging proximal algorithm will eventually be sandwiched between two extremes families of subspaces controlled by the pair $(x^\star, u^\star)$.
This identification can be exploited within our adaptive algorithm \adaalgo for solving Problem~\eqref{eq:main_problem}. Indeed, assuming that the two extreme subspaces of \eqref{eq:identification_result} coincide, the theorem says that the structure of the iterate $\S(x^k)$ will be the same as the one of the solution $\S(x^\star)$. In this case, if we choose the adaptation strategy of our adaptive algorithm \adaalgo deterministically from $\S(x^k)$, then, after a finite time with probability one,  the selection will not be adapted anymore. This allows us to recover the rate of the non-adaptive case (Theorem~\ref{th:conv_nondis}), as formalized in the next theorem.

\begin{theorem}[Improved asymptotic rate]\label{th:rate_identif}
Under the same assumptions as in Theorems~\ref{th:conv_nondis_arbitrary} and~\ref{th:aggressive}, if the solution $x^\star$ of~\eqref{eq:main_problem} verifies the qualification constraint\footnote{The qualifying constraint \eqref{eq:qualif} may seem hard to verify at first glance but for most structure-enhancing regularizers, it simplifies greatly and {reduces} to usual nondegeneracy {assumptions}. Broadly speaking, this condition simply means that the point $u^\star = x^\star - \gamma\nabla f(x^\star)$ is not \emph{borderline} to be put to an identified value by the proximity operator of the regularizer $\prox_{\gamma g}$. 
For example, {when} $g(x) = \lambda_1 \|x\|_1$, the qualifying constraint \eqref{eq:qualif} simply rewrites $x_i^\star = 0 \Leftrightarrow \nabla_i f(x^\star) \in ]-\lambda_1,\lambda_1[$; 
for $g$ is the TV-regularization \eqref{eq:TV}, the qualifying constraint means that there is no point $u$ (in any ball) around $x^\star - \gamma \nabla f(x^\star)$ such that $\prox_{\gamma g}(u)$ has a jump that $x^\star$ does not have. 
In general, this corresponds to the relative interior assumption of \cite{lewis2002active}; see the extensive discussion 
of \cite{vaiter-model-linear15}.}
\begin{align}
    \label{eq:qualif}\tag{QC}
     \S(x^\star) ~=\!\! \bigcup_{u\in\mathcal{B}(x^\star-\gamma \nabla f(x^\star),\varepsilon)} \!\S(\prox_{\gamma g}(u))
\end{align}
for any $\varepsilon>0$ small enough, then, using an adaptation deterministically computed from $(\S(x^k))$, we have \begin{align*}
     \EE[\|x^k-x^\star\|_2^2] = \mathcal{O}\left( \left( 1- \lambda_{\min} (\bP^\star) \frac{2\gamma \mu L}{\mu + L}  \right)^k \right)
 \end{align*}
 where $ \bP^\star$ is the average projection matrix of the selection 
 associated with $\S(x^\star)$.
\end{theorem}

\proof{Proof.}
Let $u^\star = x^\star-\gamma \nabla f(x^\star)$ and observe from the optimality conditions of \eqref{eq:main_problem} that $x^\star =  \prox_{\gamma g}(u^\star)$. We apply Theorem~\ref{th:identification} and the qualification condition \eqref{eq:qualif} yields that $\S(x^k)$ will exactly reach $\S(x^\star)$ in finite time. Now we go back to the proof of Theorem~\ref{th:aggressive} to see that the random variable defined by
\begin{align*}
    X^k = \left\{ \begin{array}{cl}
            x^{k_{\ell}} & \textrm{if } k\in(k_{\ell},k_{\ell}+\mathbf{c}_\ell]   \\
        x^k & \textrm{if } k\in(k_{\ell}+\mathbf{c}_\ell , k_{\ell+1}]  
    \end{array}\right. \textrm{ for some } \ell
\end{align*}
also converges almost surely to $x^\star$. Intuitively, this sequence is a replica of $(x^k)$ except that it stays fixed at the beginning of adaptation cycles when no adaptation is admitted. 
This means that $\S(X^k)$ which can be used for adapting the selection will exactly reach $\S(x^\star)$ in finite time.
From that point on, since we use an adaptation technique that deterministically relies on $\S(x^k)$, there are no more adaptations and thus the rate matches the non-adaptive one of Theorem~\ref{th:conv_nondis}.
 \endproof

This theorem means that if $g$, $\M$, and $\C$ are chosen in agreement, the adaptive algorithm \adaalgo eventually reaches a linear rate in terms of iterations as the non-adaptive \algo. In addition, the term $\lambda_{\min}(\bP)$ present in the rate now depends on the \emph{final} selection and thus on the optimal structure which is much better than the structure-agnostic selection of \algo in Theorem~\ref{th:conv_nondis}.  In the next section, we develop  practical rules for an efficient interlacing of $g$, $\M$, and $\C$.

\subsection{Identification-based subspace descent}\label{sec:adapt}

In this section, we provide practical rules to sample efficiently subspaces according to the structure identified by the iterates of our proximal algorithm. According to Theorem~\ref{th:rate_identif}, we need to properly choose $\C$ with respect to $g$ and $\M$ to have a good asymptotic regime. According to Theorem~\ref{th:aggressive},  we also need to follow specific interlacing constraints to have a good behavior along the convergence. These two aspects are discussed in Section \ref{sec:howto} and Section~\ref{sec:ex_ada}, respectively.

\subsubsection{How to update the selection}
\label{sec:howto}

We provide here general rules to sample in the family of subspaces $\C$ according to the structure identified with the family of $\M$.
To this end, we need to consider the two families $\C$ and $\M$ that closely related. We introduce the notion of generalized {complemented subspaces}.

\begin{definition}[Generalized {complemented subspaces}]
Two families of subspaces $\M = \{ \M_1,\ldots,\M_m\}$ and 
$\C = \{ \C_1,\ldots,\C_m\}$ are said to be (generalized)  {complemented subspaces} if for all $i=1,\ldots,m$
\[
  \left\{ \begin{array}{l}
    \left( \C_i \bigcap  \NC_i  \right) \subseteq  \bigcap_j \C_j  \\
      \C_i  +  \NC_i = \mathbb{R}^n
  \end{array} \right. 
\]
\end{definition}

\begin{example}[{Complemented subspaces and sparsity vectors} for axes and jumps]\label{ex:supp}
For the axes subspace set (see Section~\ref{sec:coordproj})
\begin{equation}\label{eq:Ccoord}
    \C = \{\C_1,\ldots,\C_n\} \qquad\text{ with } \C_i = \{x\in\mathbb{R}^n  : x_j = 0 ~~ \forall j \neq i \},
\end{equation}
a {complemented} identification set is 
\begin{equation}\label{eq:Mcoord}
    \NC = \{\NC_1,\ldots,\NC_n\} \quad\text{ with } \NC_i = \{x\in\mathbb{R}^n  : x_i = 0 \}, 
\end{equation} 
as $\NC_i \bigcap \C_i = \{0\} = \bigcap_j \C_j$ and $\C_i+\NC_i = \mathbb{R}^n$. 
In this case, the sparsity vector $\SC(x)$ corresponds to the \emph{support} of $x$ (indeed $[\SC(x)]_i = 0$ iff $x \in \NC_i \Leftrightarrow x_i=0$). Recall that the support of a point $x\in\mathbb{R}^n$ is defined as the size-$n$ vector $\mathrm{supp}(x)$ such that $\mathrm{supp}(x)_i = 1$ if $x_i\neq 0$ and $0$ otherwise. By a slight abuse of notation, we denote by $|\mathrm{supp}(x)|$ the size of the support of $x$, i.e. its number of non-null coordinates and $|\mathrm{null}(x)| = n - |\mathrm{supp}(x)|$.

For the jumps subspace sets (see Section~\ref{sec:var})
\begin{equation}\label{eq:Cjump}
    \C = \{\C_1,..,\C_{n-1}\} \qquad\text{ with } \C_i = \left\{x\in\mathbb{R}^n  : 
x_j = x_{j+1} \text{ for all $j\neq i$} \right\}
\end{equation}  
a  {complemented} identification set is 
\begin{equation}\label{eq:Mjump}
\NC = \{\NC_1,..,\NC_{n-1} \} \qquad\text{ with } \NC_i = \left\{x\in\mathbb{R}^n  :  x_{i} = x_{i-1} \right\},    
\end{equation}
as $\NC_i \bigcap \C_i = \mathrm{span}(\{1\}) = \bigcap_j \C_j$ and $\C_i+\NC_i = \mathbb{R}^n$. Here $\SC(x^{k})$ corresponds to the \emph{jumps} of $x$ (indeed $[\SC(x^k)]_i = 0$ iff $x^k \in \NC_i \Leftrightarrow x_i^k=x_{i+1}^k$). . The jumps of a point $x\in\mathbb{R}^n$ is defined as the vector $\mathrm{jump}(x)\in \mathbb{R}^{(n-1)}$ such that for all $i$ we have: $\mathrm{jump}(x)_i = 1$ if $x_i\neq x_{i+1}$ and $0$ otherwise.  
\end{example}

\smallskip

The practical reasoning with using {complemented} families is the following. If the subspace $\NC_i$ is identified at time $K$ (i.e. $[\S(x^k)]_i=0 \Leftrightarrow x^k\in\NC_i$ for all $k\geq K$), then it is no use to update the iterates in $\C_i$ in preference, and the next selection $\Sel_k$ should not include $\C_i$ anymore.
Unfortunately, the moment after which a subspace is definitively identified is unknown in general; however, subspaces $\NC_i$ usually show a certain stability and thus $\C_i$ may be ``less included'' in the selection. This is the intuition behind our adaptive subspace descent algorithm: when the selection $\Sel^k$ is adapted to the subspaces in $\NC$ to which $x^k$ belongs, this gives birth to an automatically adaptive subspace descent algorithm, from the generic \adaalgo.

Table~\ref{tab:comp} summarizes the common points and differences between the adaptive and non-adaptive subspace descent methods. Note that the two options introduced in this table are examples on how to generate reasonably performing admissible selections. Their difference lies in the fact that for Option\;1, the \emph{probability} of sampling a subspace outside the support is controled, while for Option\;2, the \emph{number} of subspaces is controlled (this makes every iteration computationally similar which can be interesting in practice). 
Option 2 will be 
discussed in Section\;\ref{sec:ex_ada} and 
illustrated numerically in Section~\ref{sec:num}. 

\begin{table}[h!]
    \centering
    \begin{tabular}{rl|c|c|}
      &  & (non-adaptive) subspace descent & adaptive subspace descent   \\
      &   &  \algo & \adaalgo \\
         \hline
            \multicolumn{2}{ c| }{Subspace family}  &  \multicolumn{2}{ c| }{
$ \C = \{\C_1,..,\C_c\}$} \\   
          \hline
   \multicolumn{2}{ c| }{Algorithm}  &  \multicolumn{2}{ c| }{
$ \left\{
\begin{array}{rl}
    y^k = &  \bQ\left(x^k - \gamma\nabla f\left(x^k\right)\right) \\
     z^{k} = & P_{\Sel^k} \left(y^k\right) + (I- P_{\Sel^k} ) \left(z^{k-1}\right)  \\
     x^{k+1} = & \prox_{\gamma g} \left(\bQ^{-1}\left(z^{k}\right)\right)
\end{array}\right.   $ } \\   
          \hline
     \multirow{6}{*}{Selection} &    \multirow{2}{*}{Option 1}  &    &  $\C_i\in\Sel^k$ with probability  \\
     & &  $\C_i\in\Sel^k$ with probability $p$ & $ \left\{
\begin{array}{rl}
    p &  \text{ if } x^k\in\NC_i \Leftrightarrow [\SC(x^k)]_i = 0 \\
    1 & \text{ elsewhere } 
\end{array}\right.$ \\
\cline{2-4}
 &    \multirow{4}{*}{Option 2}  &  &  Sample $s$ elements uniformly in \\
      &  &  Sample $s$ elements uniformly in $\C$ &  $\{\C_i  :  x^k\in\NC_i  \text{ i.e. }  [\SC(x^k)]_i = 0  \}$  \\
     & & & and add \emph{all} elements in \\
          & & &  $\{\C_j : x^k\notin\NC_j \text{ i.e. }  [\SC(x^k)]_j = 1 \}$ \\
     \hline 
    \end{tabular}
    \caption{Strategies for non-adaptive vs.\;adaptive algorithms}
    \label{tab:comp}
\end{table}

Notice that, contrary to the importance-like adaptive algorithms of \cite{stich2017safe} for instance, the purpose of these methods is not to adapt each subspace probability to local \emph{steepness} but rather to adapt them to the current \emph{structure}. This is notably due to the fact that local steepness-adapted probabilities can be difficult to evaluate numerically and that in heavily structured problems, adapting to an ultimately very sparse structure already reduces drastically the number of explored dimensions, as suggested in \cite{grishchenko2018asynchronous} for the case of coordinate-wise projections.

\subsubsection{Practical examples and discussion}\label{sec:ex_ada}

We discuss further the families of subspaces of Example\;\ref{ex:supp} when selected with Option 2 of Table~\ref{tab:comp}.

\paragraph{Coordinate-wise projections}\label{Pex:l1}

Using the subspaces \eqref{eq:Ccoord} and \eqref{eq:Mcoord}, a practical adaptative coordinate descent can be obtained from the following reasoning at each adaptation time $k = k_{\ell-1}$:
\begin{itemize}
    \item Observe $\SC(x^{k})$ i.e. the support of $x^{k}$.
    \item Take all coordinates in the support and randomly select $s$ coordinates outside the support. 
    Compute\footnote{Let us give a simple example in $\mathbb{R}^4$:
    \begin{align*}
        \text{for } x^k = \left( \begin{array}{c}
            1.23 \\ -0.6 \\ 0 \\ 0 
        \end{array}  \right) \text{,  } \SC(x^{k}) =  \left( \begin{array}{c}
            1 \\ 1 \\ 0 \\ 0 
        \end{array}  \right) \text{, ~then~ } \begin{array}{l}
            \Prob[\C_1 \subseteq \Sel^k ] = \Prob[\C_2 \subseteq \Sel^k ] = 1 \\
            \Prob[\C_3 \subseteq \Sel^k ] = \Prob[\C_4 \subseteq \Sel^k ] = p_\ell :=  s/|\mathrm{null}(x^k)|= s/2
        \end{array}   
    \end{align*}  
    \begin{align*}
\bP_\ell = \left( \begin{array}{cccc}
            1 & & & \\  & 1 & &  \\   & & p_\ell &  \\   & & & p_\ell
        \end{array} \right) ~~~~ \bQ_\ell = \left( \begin{array}{cccc}
            1 & & & \\  & 1 & &  \\   & & 1/\sqrt{p_\ell} &  \\   & & & 1/\sqrt{p_\ell}
        \end{array} \right) ~~~~  \bQ_\ell^{-1} = \left( \begin{array}{cccc}
            1 & & & \\  & 1 & &  \\   & & \sqrt{p_\ell} &  \\   & & & \sqrt{p_\ell}
        \end{array} \right)~~~~~
    \end{align*}} associated $\bP_\ell  $, $\bQ_\ell$, and $\bQ_\ell^{-1}$. Notice that $ \lambda_{\min}(\bP_{\ell}) = p_{\ell} = s/|\mathrm{null}(x^k)|$.
    \item  Following the rules of Theorem~\ref{th:aggressive}, compute
    \begin{align*}
        \mathbf{c}_\ell = \left\lceil \frac{\log\left(\|  \bQ_\ell   \bQ_{\ell-1}^{-1} \|_2^2 \right) + \log(1/(1-\beta)) }{\log( 1/(1-\alpha_{\ell-1}))}\right\rceil \qquad
        \text{with~}\alpha_{\ell-1} =2 p_{\ell-1} \gamma\mu L /(\mu+L)
    \end{align*}  
    for some small fixed $0<\beta \leq 2 \gamma\mu L /(n(\mu+L)) \leq\inf_{\ell} \alpha_\ell $.

    Apply the new sampling after $\mathbf{c}_\ell$ iterations (i.e. $k_{\ell} = k_{\ell-1}+\mathbf{c}_\ell$).
\end{itemize}

\smallskip
Finally, we notice that the above strategy with Option\;2 of Table~\ref{tab:comp} produces moderate adaptations as long as the iterates are rather dense. To see this, observe first that $\bQ_\ell   \bQ_{\ell-1}^{-1}$ is a diagonal matrix, the entries of which depend on the support of the corresponding coordinates at times ${k_{\ell -1}}$ and ${k_{\ell -2}}$. More precisely, the diagonal entries are described in the following table:

\begin{center}
\smallskip
\begin{tabular}{p{0.1\textwidth}|p{0.1\textwidth}|c}
   \multicolumn{2}{c}{$i$ is in the support at} & \\
   ${k_{\ell -1}}$ & ${k_{\ell -2}}$ &  $\left[\bQ_\ell   \bQ_{\ell-1}^{-1}\right]_{ii}$\\
    \hline
   yes & yes &  $1$  \\
    \hline
   no & yes &  $\frac{1}{p_{\ell}} = \frac{|\mathrm{null}(x^{k_{\ell-1}})|}{s}$ \\
    \hline
    yes  & no & ${p_{\ell-1}} = \frac{s}{|\mathrm{null}(x^{k_{\ell-2}})|}$\\
    \hline
    no & no & $\frac{p_{\ell-1}}{p_{\ell}} = \frac{|\mathrm{null}(x^{k_{\ell-1}})|}{|\mathrm{null}(x^{k_{\ell-2}})|}$
\end{tabular}
\smallskip
\end{center}
Thus, as long as the iterates are not sparse (i.e. in the first iterations, when $|\mathrm{null}(x^k)|\approx s$ is small), the adaptation cost is moderate so the first adaptations can be done rather frequently. Also, in the frequently-observed case when the support only decreases ($\S(x^{k_{\ell-2}})\leq \S(x^{k_{\ell-1}})$), the second line of the table is not active and thus $\|\bQ_\ell  \bQ_{\ell-1}^{-1}\| = 1$, so the adaptation can be done without waiting.

\paragraph{Vectors of fixed variations}
\label{Pex:TV}

The same reasoning as above can be done for vectors of fixed variation by using the families \eqref{eq:Cjump} and \eqref{eq:Mjump}. At each adaptation time $k = k_{\ell-1}$:
\begin{itemize}
    \item Observe $\SC(x^{k})$ i.e. the \emph{jumps} of $x$;
    \item The adapted selection consists in selecting all jumps present in $x^k$ and randomly selecting $s$ jumps that are not in $x^k$. Compute $\bP_\ell  $, $\bQ_\ell$, and $\bQ_\ell^{-1}$ (to the difference of coordinate sparsity they have to be computed numerically). 
    \item For a fixed $\beta>0$, compute
    \begin{align*}
        \mathbf{c}_\ell = \left\lceil \frac{\log\left(\|  \bQ_\ell   \bQ_{\ell-1}^{-1} \|_2^2 \right) + \log(1/(1-\beta)) }{\log( 1/(1-\alpha_{\ell-1}))}\right\rceil.
    \end{align*}  
     Apply the new sampling 
     after $\mathbf{c}_\ell$ iterations (i.e. $k_{\ell} = k_{\ell-1}+\mathbf{c}_\ell$).
\end{itemize}

\section{Numerical illustrations}\label{sec:num}

We report preliminary numerical experiments illustrating the behavior of our randomized proximal algorithms on standard problems involving $\ell_1$/TV regularizations. We provide an empirical comparison of our algorithms with the standard proximal (full and coordinate) gradient algorithms and a recent proximal sketching algorithm.

\subsection{Experimental setup}

We consider the standard regularized logistic regression with three different regularization terms, which can be written for given $(a_i,b_i)\in\RR^{n+1}$ ($i=1,\ldots,m$) and parameters $\lambda_1,\lambda_2>0$
\begin{subequations}
\begin{align}
 &+ ~\lambda_1\!\left\|x\right\|_1\label{eq:logl1}\\
 \min_{x\in \RR^n}~~~\frac{1}{m}\sum\limits_{i=1}^m \log\left(1 + \exp\left(-b_ia_i^\top x\right)\right) + \frac{\lambda_2}{2}\|x\|_2^2 &+ ~  \lambda_1\!\left\|x\right\|_{1,2}\label{eq:logl12}\\
 &+~\lambda_1\!\mathbf{TV}(x)\label{eq:logtv}
\end{align}
\end{subequations}
We use two standard data-sets from the LibSVM repository: the \emph{a1a} data-set ($m=1,605$ $n=123$) for the $\mathbf{TV}$ regularizer, and the \emph{rcv1\_train} data-set ($m=20,242$ $n=47,236$) for the $\ell_1$ and $\ell_{1,2}$ regularizers. We fix the parameters $\lambda_2=1/m$ and $\lambda_1$ to reach a final sparsity of roughly 90\%. 

The subspace collections are taken naturally adapted to the regularizers: by coordinate for \eqref{eq:logl1} and \eqref{eq:logl12}, and by variation for \eqref{eq:logtv}. The adaptation strategies are the ones described in Section~\ref{sec:ex_ada}.

We consider five algorithms:\\[0.2cm]
\resizebox{\linewidth}{!}{
\begin{tabular}{|c|c|c|c|}
\hline
    Name & Reference & Description & Randomness  \\
    \hline
   \pgd  &  &  vanilla proximal gradient descent  & None \\
   x\footnotemark \, \rpcd & \cite{nesterov2012efficiency} & standard proximal coordinate descent &  x coordinates selected for each update \\
   x  \sega & \cite{hanzely2018sega} & Algorithm \sega~with coordinate sketches &  $\rank(S^k) = \text{x}$ \\
   x  \algo &  Algorithm \ref{alg:strata_nondis} & (non-adaptive) random subspace descent & Option 2 of Table~\ref{tab:comp} with $s =  \text{x}$ \\
   x  \adaalgo & Algorithm \ref{alg:ada_strata_nondis} & adaptive random subspace descent &  Option 2 of  Table~\ref{tab:comp} with $s =  \text{x}$\\
 \hline
\end{tabular}}
\footnotetext{In the following, x is often given in percentage of the possible subspaces, i.e. x\% of $|\mathcal{C}|$, that is x\% of $n$ for coordinate projections and x\% of $n-1$ for variation projections.}

\vspace*{0.5cm}

For the produced iterates, we measure the sparsity of a point $x$ by $\|\S(x_k)\|_1$, which corresponds to the size of the supports for the $\ell_1$ case and the number of jumps for the TV case. We also consider the quantity: 
\[
\text{Number of subspaces explored at time $k$} \displaystyle ~=~ \sum^k_{t=1} \|\S(x^t)\|_1.
\]
We then compare the performance of the algorithms on three criteria:
\begin{itemize}
    \item functional suboptimality vs iterations (standard comparison);
    \item size of the sparsity pattern vs iterations (showing the identification properties);
    \item functional suboptimality vs number of subspaces explored (showing the gain of adaptivity).
\end{itemize}

\subsection{Illustrations for coordinate-structured problems}

\subsubsection{Comparison with standard methods}

We consider first $\ell_1$-regularized logistic regression \eqref{eq:logl1}; in this setup, the non-adaptive \algo~boils down to the usual randomized proximal gradient descent (see Section\;\ref{sec:coordproj}). We compare the proximal gradient to its adaptive and non-adaptive randomized counterparts.

First, we observe that the iterates of \pgd and \adaalgo coincide. This is due to the fact that the sparsity of iterates only decreases ($\S(x_k)\leq \S(x_{k+1})$) along the convergence, and according to Option 2 all the non-zero coordinates are selected at each iteration and thus set to the same value as with \pgd. However, a single iteration of $10\%$-\adaalgo costs less in terms of number of subspaces explored, leading the speed-up of the right-most plot.  Contrary to the adaptive \adaalgo, the structure-blind \algo~identifies much later then \pgd and shows poor convergence.

\begin{figure}[H]
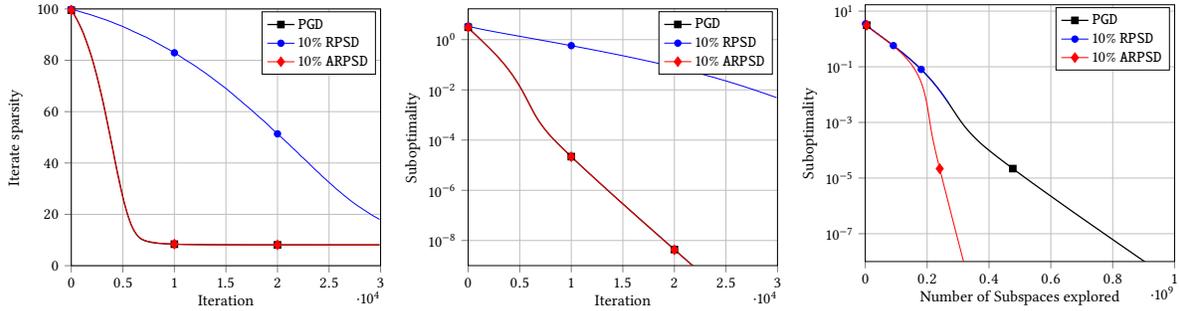

\begin{center}
\scalebox{.6}{\input{Article/rcv1_sparsity_vs_iteration_l1__algorithms_comparison.tex}}
\scalebox{.6}{\input{Article/rcv1_suboptimality_vs_iteration_l1_different_amount_of_coordinates.tex}}
\scalebox{.6}{\input{Article/rcv1_suboptimality_vs_sparsity_l1_different_amount_of_coordinates.tex}}
\end{center}
\caption{$\ell_1$-regularized logistic regression \eqref{eq:logl1}
}
\label{fig:rcv1}
\end{figure}

\subsubsection{Comparison with \sega}\label{sec:num:sega}

\begin{figure}[H]
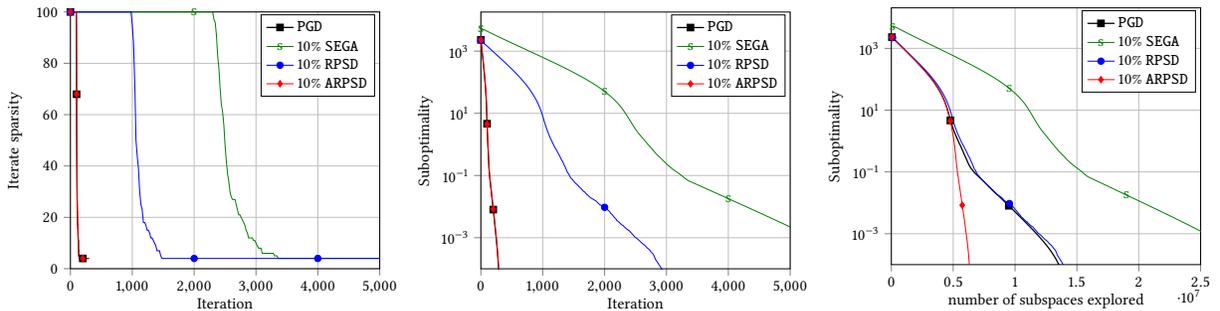

\begin{center}
\scalebox{.6}{\input{Article/sega_sparsity_vs_iteration_tv__algorithms_comparison.tex}}
\scalebox{.6}{\input{Article/sega_suboptimality_vs_iteration_tv_algorithms_comparison.tex}}
\scalebox{.6}{\input{Article/sega_suboptimality_vs_sparsity_tv__algorithms_comparison.tex}}
\end{center}
\caption{$\ell_{1,2}$ regularized logistic regression \eqref{eq:logl12}
}
\label{fig:rcv1_l12}
\end{figure}

In Figure~\ref{fig:rcv1_l12}, we compare \adaalgo~algorithm with \sega~algorithm featuring coordinate sketches\;\cite{hanzely2018sega}. While the focus of \sega~is not to produce an efficient coordinate descent method but rather to use sketched gradients, \sega~and \algo~are similar algorithmically and reach similar rates (see Section \ref{sec:comparison}). As mentioned in \cite[Apx.\;G2]{hanzely2018sega}, \sega~is slightly slower than plain randomized proximal coordinate descent (10\% \algo) but still competitive, which corresponds to our experiments. Thanks to the use of identification, \adaalgo~shows a clear improvement over other methods in terms of efficiency with respect to the number of subspaces explored.

\subsection{Illustrations for total variation regularization}

\begin{figure}[H]
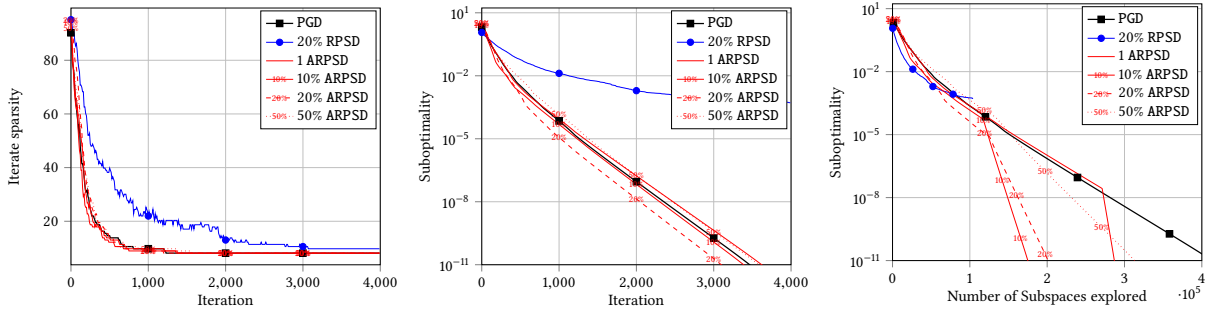

\scalebox{.6}{\input{Article/a1a_sparsity_vs_iteration_tv_different_amount_of_coordinates.tex}}
\scalebox{.6}{\input{Article/a1a_suboptimality_vs_iteration_tv_different_amount_of_coordinates.tex}}
\scalebox{.6}{\input{Article/a1a_suboptimality_vs_sparsity_tv_different_amount_of_coordinates.tex}}
\caption{1D-TV-regularized logistic regression \eqref{eq:logtv}}
\label{fig:a11}
\end{figure}

We focus here on the case of total variation \eqref{eq:logtv} which is a typical usecase for our adaptive algorithm and subspace descent in general. Figure~\ref{fig:a11} displays a comparison between the vanilla proximal gradient and various versions of our subspace descent methods.

We observe first that  \algo, not exploiting the problem structure, fails to reach satisfying performances as it identifies lately and converges slowly. In contrast, the adaptive versions \adaalgo perform similarly to the vanilla proximal gradient in terms of sparsification and suboptimality with respect to iterations. As a consequence, in terms of number of subspaces explored, \adaalgo becomes much faster once a near-optimal structure is identified. More precisely, all adaptive algorithms (except 1 \adaalgo, see the next paragraph) identify a subspace of size $\approx 8\%$ (10 jumps in the entries of the iterates) after having explored around $10^5$ subspaces. Subsequently, each iteration involves a subspace of size 22,32,62 (out of a total dimension of 123) for 10\%,20\%,50\% \adaalgo respectively, resulting in the different slopes in the red plots on the rightmost figure.

\begin{figure}[H]
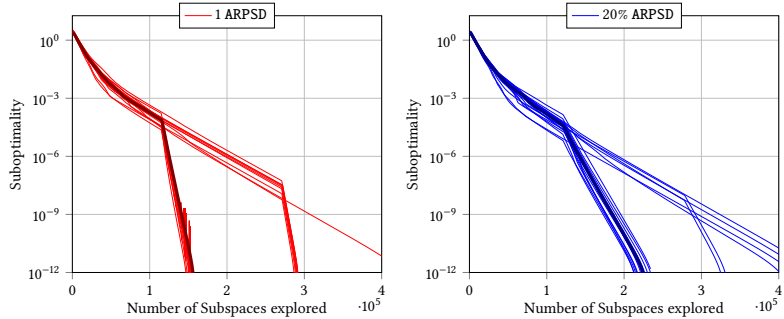

\begin{center}
\scalebox{.6}{\input{Article/a1a_expectation_1_tv_different_amount_of_coordinates.tex}}
\scalebox{.6}{\input{Article/a1a_expectation_20_tv_different_amount_of_coordinates.tex}}
\end{center}
\caption{20 runs of \adaalgo and their median (in bold) on 1D-TV-regularized logistic regression \eqref{eq:logtv}
}
\label{fig:tvcomp}
\end{figure}

Finally, Figure~\ref{fig:tvcomp} displays 20 runs of 1 and 20\% \adaalgo as well as the median of the runs in bold. We notice that more than 50\% of the time, a low-dimensional structure is quickly identified (after the third adaptation) resulting in a dramatic speed increase in terms of subspaces explored. However, this adaptation to the lower-dimensional subspace might take some more time (either because of poor identification in the first iterates or because a first heavy adaptation was made early and a pessimistic bound on the rate prevents a new adaptation in theory). Yet, one can notice that these adaptations are more stable for the 20\% than for the 1 \adaalgo, illustrating the ``speed versus stability'' tradeoff in the selection.

%
%
%

%

\appendix

\section{Convergence in the non-strongly convex case}

In this appendix, we study the convergence of the subspace descent algorithms, when the smooth function $f$ is convex but not strongly convex. Removing the strong convexity from Assumption~\ref{hyp:f}, we need existence of the optimal solutions of \eqref{eq:main_problem} and thus we make the following assumption.

\begin{assumption}\label{hyp:f2}
The function $f$ is convex $L$-smooth and the function $g$ is convex, proper, and lower-semicontinuous. Let $X^\star\neq\emptyset$ denote the set of minimizers of Problem~\eqref{eq:main_problem}.
\end{assumption}

With Assumption~\ref{hyp:f2} replacing Assumption~\ref{hyp:f},
the convergence results (Theorems\;\ref{th:conv_nondis} and\;\ref{th:conv_nondis_arbitrary}) extend from similar rationale. Let us here formalize the result and its proof for the non-adaptive case: the next theorem establishes the convergence of \algo, still with the usual fixed stepsize in $(0,2/L)$.

\begin{theorem}[\algo~convergence]\label{th:conv}
Let Assumptions \ref{hyp:f2} and \ref{hyp:main} hold. Then, for any with $\gamma\in(0,2/L)$, the sequence $(x^k)$ of the iterates of \algo converges almost surely to a point in the set $X^\star$ of the minimizers of \eqref{eq:main_problem}.
\end{theorem}

To prove this result, one can first notice that Lemma~\ref{lm:removing_exp} still holds, contrary to Lemma~\ref{lm:bub}. Thus, let us provide a replacement for Lemma~\ref{lm:bub} in the non-strongly convex setup.

\begin{lemma}\label{lm:bub2}
If Assumptions \ref{hyp:f2} and \ref{hyp:main} holds, then for $\gamma \in (0,2/L)$ and for any $x^\star \in X^\star$ (with associated $z^\star = y^\star = \bQ\left(x^\star - \gamma\nabla f\left(x^\star\right)\right)$ ), one has 
    \begin{equation*}
\|y^{k}-y^\star\|_{\bP}^2 -  \|z^{k-1}-z^\star\|_{\bP}^2 \leq - \frac{2- \gamma L}{\gamma L}  \|  \nabla f(x^k) -\nabla f(x^\star) \|_2^2.
\end{equation*}
\end{lemma}

\proof{Proof.}
Using the same arguments as in the proof of Lemma~\ref{lm:bub}, we can also show that
\begin{align}
\label{eq:A1}  \|y^{k}-y^\star\|_{\bP}^2  &= \left\| x^k - \gamma\nabla f(x^k) - ( x^\star - \gamma\nabla f(x^\star) ) \right\|_2^2 ;\\
\label{eq:A2}   \text{ and  } ~~~ \|x^k - x^\star\|_2^2 &\leq  \|z^{k-1}-z^\star\|_{\bP}^2.
\end{align}
Now, using the Baillon-Haddad theorem (see \cite[Cor. 18.16]{bauschke2011convex}), for $\gamma \in (0,2/L)$, one has 
\begin{align*}
\| x^k - \gamma\nabla f(x^k) - ( x^\star - \gamma\nabla f(x^\star) ) \|_2^2 &\leq  \| x^k  -  x^\star  \|_2^2 - \frac{2- \gamma L}{\gamma L} \|  \nabla f(x^k) -\nabla f(x^\star) \|_2^2.
\end{align*}
Combining with \eqref{eq:A1},\eqref{eq:A2} directly leads to the result.
 \endproof

\proof{Proof.}\textit{(of Theorem\;\ref{th:conv})} 
Combining Lemmas~\ref{lm:removing_exp} and \ref{lm:bub2}, we get
for any $x^\star\in X^\star$ and associated $z^\star = \bQ\left(x^\star - \gamma\nabla f\left(x^\star\right)\right)$
\begin{align}
\label{eq:basis_nonstrong}
    \EE\left[\|z^{k} - z^\star\|_2^2\,|\,\FF^{k-1}\right] \leq \|z^{k-1}-z^\star\|_2^2 - \frac{2- \gamma L}{\gamma L} \|  \nabla f(x^k) -\nabla f(x^\star) \|_2^2.
\end{align}
Taking the expectation on both sides and telescoping, we get that $\EE[\sum_{k=1}^\infty \|  \nabla f(x^k) -\nabla f(x^\star) \|_2^2] < \infty$ and thus $\nabla f(x^k) \to \nabla f(x^\star)$ with probability one. 

Eq.~\eqref{eq:basis_nonstrong} also implies that, 
as in the strongly convex case, the sequence $(\|z^{k} - z^\star\|_2^2)$ 
is a non-negative super-martingale with respect to the filtration $(\FF^k)$ and thus converges to a finite random variable (in fact, that is a common observation for randomized monotone operators; see e.g. \cite[Apx.~B]{bianchi2016coordinate}).
As a consequence, the sequence $(z^k)$ is bounded almost surely. Let $\overline{z}$
be an accumulation point of $(z^k)$; it verifies $\nabla f( \prox_{\gamma g}(\bQ^{-1} \overline{z}) ) = \nabla f(x^\star) $ and is thus in $Z^\star = \{ \bQ\left(x - \gamma\nabla f(x)\right) : x\in X^\star\}$. Denote $\overline{x} \in X^\star$ such that $\overline{z} = \bQ\left(\overline{x} - \gamma\nabla f\left(\overline{x}\right)\right)$.

Using for $\overline{x}$ the same rationale as above for $x^\star$, we can prove that the sequence $\|z^{k} - \overline{z}\|_2^2$  converges. Therefore, we deduce that with probability one, $\lim \|z^{k} - \overline{z}\|_2^2 = \lim\inf \|z^{k} - \overline{z}\|_2^2 = 0$. 
This shows that $(z^k)$ converges almost surely to $ \overline{z}$.  Applying the map $ \prox_{\gamma g} \circ \bQ^{-1}$ to this result leads to the claimed result.
 \endproof

\section*{Acknowledgments.}
The authors benefited from the support of IDEX Grenoble Alpes IRS grant \emph{DOLL}.




\bibliographystyle{ormsv080}
\bibliography{references.bib}

\end{document}